\documentclass[a4paper, reqno, 11pt]{amsart}
\usepackage[utf8]{inputenc}
\usepackage{amsmath}
\usepackage{amssymb}
\usepackage{amsthm}
\usepackage{diffcoeff}
\usepackage{bm}
\usepackage{color, xcolor}
\usepackage{graphicx}
\usepackage{mathtools}
\usepackage{appendix}
\usepackage{subcaption}
\usepackage{makecell}
\usepackage{multirow}
\usepackage{booktabs}
\usepackage{lipsum}
\usepackage{hyperref}
\hypersetup{
    colorlinks,
    linkcolor = {blue!50!black},
    citecolor = {blue!50!black},
    urlcolor  = {blue!80!black}
}

\newcommand{\bB}[1]{\boldsymbol{#1}}

\DeclareFontFamily{U}{mathx}{}
\DeclareFontShape{U}{mathx}{m}{n}{<-> mathx10}{}
\DeclareSymbolFont{mathx}{U}{mathx}{m}{n}
\DeclareMathAccent{\widehat}{0}{mathx}{"70}
\DeclareMathAccent{\widecheck}{0}{mathx}{"71}

\DeclareMathOperator{\sgn}{sgn}
\renewcommand{\div}{\mathop{\mathrm{div}}}
\newcommand{\bdiv}{\mathop{\mathbf{div}}}

{\theoremstyle{remark} \newtheorem{remark}{Remark}}

\title{A moment model of shallow granular flows with variable friction laws}

\author[J. Careaga]{Julio Careaga$^{\lowercase{\textrm{a}}}$}
\author[Q. Huang]{Qian Huang$^{\lowercase{\textrm{b}}}$}
\author[J. Koellermeier]{Julian Koellermeier$^{\lowercase{\textrm{a,c}}}$}

\date{\today}

\makeatletter
\def\@settitle{\begin{center}%
  \normalfont\LARGE\bfseries
  \@title
  \end{center}}
\makeatother

\thanks{\noindent$^{\text{a}}$Bernoulli Institute, University of Groningen, Nijenborgh 9, 9747 AG Groningen, The Netherlands.}
\thanks{$^{\text{b}}$Institute of Applied Analysis and Numerical Simulation, University of Stuttgart, 70569 Stuttgart, Germany.}
\thanks{$^{\text{c}}$Department of Mathematics, Computer Science and Statistics, Ghent University, 9000 Ghent, Belgium.}
\thanks{\text{Email addresses:}
{\tt j.c.careaga.solis@rug.nl},
{\tt qian.huang@mathematik.uni-stuttgart.de},
{\tt j.koellermeier@rug.nl}}

\usepackage{titlesec}

\titleformat{\title}[hang]{\it}{}{0.7cm}{}
\titleformat{\section}[hang]{\Large\bfseries}{\thesection}{0.7cm}{}
\titlespacing*{\section}{0cm}{0.7cm}{0.5cm}
\titleformat{\subsection}[hang]{\bfseries}{\thesubsection}{0.2cm}{} %%filcenter
\titleformat{\subsubsection}[hang]{\bfseries}{\thesubsubsection}{0.2cm}{}

\begin{document}

\maketitle

\begin{abstract}
In this work, we develop a modelling framework for granular flows based on the shallow water moment equations on inclined planes.
Under the assumption of a polynomial expansion of the velocity field, the model extends the classical shallow water equations to vertically variable velocity profiles. The friction effects, which are captured through the strain-rate tensor, are incorporated into the model in two terms, the bulk and bottom friction. We propose a modelling procedure to incorporate general friction laws into our framework and exemplify this combining the Manning, Coulomb, Savage-Hutter, and $\mu(I)$-rheology friction models in our modeling framework. Moreover, we develop a path-conservative finite volume numerical scheme based on the polynomial viscosity matrix method to properly handle the stiffness of the source terms.
Numerical simulations are presented for different models of friction, including the case of wet-dry fronts.

\bigskip

\smallskip

\end{abstract}

\section{Introduction}

Fluid and granular flows, under the assumption of shallowness, are widely modelled by the so-called shallow water equations (SWE). Although in both flow regimes physical properties such as viscosity and density change, the SWE allow to describe both the water height (or layer thickness) and the depth-averaged velocity profile. In the case of flows over inclined planes, gravity is the main driving force, while friction effects act against the flow.
The modelling of such flows is essential in the understanding of natural disasters such as avalanches \cite{Fei2020,Kelfoun2005,Mangeney2003,Savage1991}, extrusions of lava flows (pyroclastic flow) \cite{Neglia2021}
or river floodings \cite{Bates2022,Kitsikoudis2020}. Granular flows are also common in industrial applications, where grains of different types are discharged from a hopper onto an inclined conveyor \cite{Cleary2017}.

Even though the depth-averaged approach of the SWE plays an important role in reducing the complexity of the underlying Navier-Stokes equations, the simple SWE are of limited use for studying velocity fields varying with respect to depth. However, experiments have demonstrated that the velocity profiles in shallow water flows tend to vary significantly depending on the vertical coordinate, see \cite{Kern2010,Sanvitale2016}. An extension of the SWE that overcomes this issue is given by the shallow water moment equations (SWME) \cite{kowalski2018moment}, in which the velocity field is assumed to be polynomial. The coefficients of the polynomial expansion are called moments, which also become additional variables. A closed system of equations is obtained by Galerkin projection onto test functions, leading to the resulting PDE system. In this model, depending on the number of moments (i.e.~the highest polynomial degree), the velocity profiles can include not only the classical linear and parabolic profiles, but more involved polynomial expressions for a high number of moments leading to increased accuracy in numerical simulations \cite{kowalski2018moment}. Several variants of the SWME have been proposed in the literature, all focusing on different treatments of the transport operator: The hyperbolic shallow water moment equations (HSWME) \cite{Koellermeier2020} which are a modified version of the SWME such that the PDE system is guaranteed hyperbolic; The shallow water linearised moment equations \cite{Pimentel2020}, which allow for a simpler analysis of the steady states of the system; A non-hydrostatic pressure version is proposed in \cite{Scholz2023}; 2D models are analysed in \cite{bauerle2024rotationalinvariancehyperbolicityshallow,Verbiest2025}.

Alternative extensions to SWEs are the multilayer shallow water models \cite{FernandezNieto2016,FernandezNieto2018,Garres-Diaz2021}, which have been used to simulate tsunamis and coastal forest interaction \cite{Burger2025a}, sediment transport \cite{Burger2025b}, and reactive sedimentation processes \cite{Careaga2024}. These multilayer models consider discontinuous velocity profiles between layers, which requires establishing sophisticated jump conditions between layers. The two approaches (multilayer and moment models) are not mutually restrictive, and their possible combinations are also of interest \cite{Garres2023}.

Friction effects are crucial in granular flows and caused by the interaction between the flow and the bottom surface as well as viscous stresses inside the fluid. These effects are typically incorporated to the SWME as a source term of the system. The standard friction law, also included in the original SWME \cite{kowalski2018moment}, is obtained as a result of considering a Newtonian flow and simple slip boundary conditions at the bottom surface. For more general rheologies, the viscous stress tensor can be given by more complicated expressions in terms of the unknowns, leading to nonlinear friction terms that need special treatment in numerical simulations. In addition, depending on the type of surface in which the flow is moving, the bottom friction may vary and therefore different friction laws have to be considered in the moment models. So far, a treatment of general friction terms is lacking for the SWME models.

Different friction laws are described in the literature. The so-called Savage-Hutter friction model \cite{Savage1989,Savage1991} relates the ratio between the shear stress and normal stress through a friction angle. Furthermore, the authors assume a cohesionless Mohr--Coulomb type material and introduce a constant internal friction angle to model the shear stress.
The Coulomb friction \cite{Coulomb1776} states that the shear stress is proportional to the viscosity and pressure.
An empirical formula is the so-called Manning friction \cite{Manning1891}, which has already been implemented in the context of a particular SWME model in \cite{Garres2020}. A well established friction law for granular flows employs the so-called $\mu(I)$-rheology \cite{Jop2005,Jop2006}, where the viscosity depends on the shear rate.

In this work, we develop a modelling framework for granular flows based on the SWME on inclined planes including bottom and bulk friction terms and variable topography. The goal is to obtain a general framework incorporating both complex friction and vertical variations to achieve accurate results in simulating shallow granular flows. This new modelling framework, while general, allows the explicit analytical derivation of model equations. We exemplify our modelling procedure with various friction terms proposed in the literature, with special emphasis on the $\mu(I)$-rheology, for which we study equilibrium states in the case of linear velocity profiles. In addition, we propose a path-conservative numerical scheme tailored to handle the distinctive non-conservative products arising in the moment models and the nonlinearities and integrals in the friction (source) terms. In addition, we address the case of wet-dry fronts that are characteristic in the simulation of avalanches in granular flows.

The paper is organised as follows. In Section~\ref{sec:SWmodels}, we derive the SWME for granular flows on inclined planes. The friction terms are then obtained as a result of treating the viscous stress tensor. The modelling procedure for handling the friction terms is presented in Section~\ref{sec:framework}, and this procedure is applied to various friction laws available in the literature. In Section~\ref{sec:granular}, we derive the first SWME model for granular flows with $\mu(I)$-rheologies. The numerical scheme is proposed in Section~\ref{sec:num:scheme} and numerical simulations are shown in Section~\ref{sec:simulations}. Finally, the conclusions are included in Section~\ref{sec:conclusions}.

\section{A general friction shallow water moment model} \label{sec:SWmodels}

In this section, we use appropriate scalings and a polynomial expansion of the velocity to derive the two-dimensional shallow water moment models for granular flows incorporating a general friction term, which corresponds to a source term depending on the unknowns of the system.  Throughout this section and in the rest of the paper, we restrict our work to the case of one horizontal and one vertical coordinate, and note that the case with two horizontal spatial coordinates is straightforward. Different specific friction models will be discussed in the following sections.

\subsection{Two-dimensional incompressible Navier Stokes Equations}
\begin{figure}[t!]
	\centering
    \includegraphics[width=0.6\textwidth]{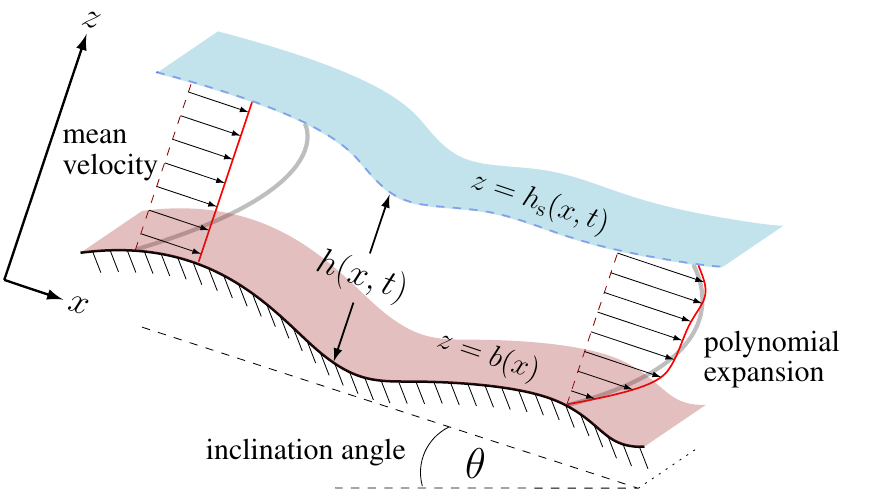}
	\caption{Shallow flow geometry on an inclined plane with angle $\theta$ and bottom topography $b(x)$ and two velocity profiles: constant (left) and more accurate polynomial (right). \label{fig:Schematic}}
\end{figure}

We study granular flows over inclined beds using shallow flow models and different rheologies. We begin by assuming that the bottom topography is described by a curve $z=b(x)$, where the pair $(x,z)$ corresponds to the horizontal and vertical coordinates, respectively, defined on an inclined plane with fixed angle $\theta\in[0,\tfrac{\pi}{2})$, see Figure~\ref{fig:Schematic}. 
We define the two-dimensional fluid velocity $\bB{u} = (u,w)^{\rm T}\in \mathbb{R}^2$ and pressure $p\in \mathbb{R}$, both functions of the spatial coordinates $(x,z)$ and time $t\geq 0$. Then, the two-dimensional incompressible Navier-Stokes equations that determine the evolution of the pair $(\bB{u}, p)$ are given by
\begin{subequations} \label{eq:govdim}
\begin{align}
\div(\bB{u}) &=0,  \label{eq:govdim:mass} \\[1ex]
  \diffp{\bB{u}}{t} + \bB{u}\cdot \nabla \bB{u} &= g \begin{pmatrix}\sin(\theta)\\ -\cos(\theta)\end{pmatrix} - \frac{1}{\rho} \nabla p + \frac{1}{\rho}\bdiv(\bB{\tau}), \label{eq:govdim:mom}
\end{align}
\end{subequations}
where $\nabla(\cdot)$ and $\div(\cdot)$ are the gradient and divergence operator for vectors, respectively, and $\bdiv(\cdot)$ is the divergence operator for tensors applying row-wise.
The free surface is denoted by $z=b(x)+h(x,t)= h_s(x,t)$, with the water height $h = h(x,t)$. The deviatoric viscous stress tensor $\bm{\tau} \in \mathbb{R}^{2\times 2}$ is defined as
\begin{align} \label{eq:tau:components}
\bm{\tau} = 
\begin{pmatrix}
\tau_{xx} & \tau_{xz} \\
\tau_{zx} & \tau_{zz}
\end{pmatrix},
\end{align}
and the total stress tensor reads $\bm \sigma = \bm \tau -p \bm I$, where $\bm I$ is the identity matrix in $\mathbb{R}^{2\times 2}$.
The stress tensor is determined by the symmetric strain-rate tensor $\bm{D}(\bm u)= \nabla \bB{u} + (\nabla\bB{u})^{\rm{T}}$ and the viscosity $\eta$ using
\begin{align}
   \bm \tau = \eta \bm D(\bm u). \label{eq:constit}
\end{align}
The precise definition of the viscosity $\eta$ will depend on the specific rheology used, which can be Newtonian flow or granular flow, among others. See Section~\ref{sec:framework} for further details.

System~\eqref{eq:govdim} needs to be supplemented with boundary conditions in the vertical and horizontal directions. The variations at the free-surface $z=h_{\rm s}$ and bottom $z=b$ lead to the following kinematic boundary conditions:
\begin{subequations} \label{eq:kinematic_BQ}
\begin{align}
  \diffp{h_{\rm s}}{t} + u(x,h_{\rm s},t) \diffp{h_{\rm s}}{x} &= w(x,h_{\rm s},t), \label{eq:kinematic:hs}\\
  u(x,b,t) \diffp{b}{x} &= w(x,b,t). \label{eq:kinematic:b}
\end{align}
\end{subequations}
In addition, dynamic boundary conditions include the stress at the surface and bottom of the fluid. To introduce these boundary conditions, we first define the downward normal and tangential vector at $z = b$, and the upward normal and tangential vector at $z=h_{\rm s}$, respectively by
\begin{alignat}{3}
  \bB{n}_{\rm b} &= \big\|(\partial_x b (x),-1)\big\|^{-1} \big(\partial_x b (x),-1\big)^{\rm T},\qquad &
  \bB{t}_{\rm b} &= \|(1, \partial_x b (x))\|^{-1}\big(1,\partial_x b (x)\big)^{\rm T}, & \label{eq:normal:tangent:b}\\
  \bB{n}_{\rm s} &= \|(-\partial_x h_{\rm s}(x),1)\|^{-1}\big(-\partial_x h_{\rm s}(x),1\big)^{\rm T}, \qquad &
  \bB{t}_{\rm s} &= \|(1, \partial_x h_{\rm s}(x))\|^{-1}\big(1,\partial_x h_{\rm s}(x)\big)^{\rm T}.& \label{eq:normal:tangent:s}
  \end{alignat}
We then impose the following dynamic boundary conditions for the stress at the bottom $z=b$ and surface $z=h_{\rm s}$:
\begin{subequations}
\begin{align}
(\bm \sigma \bm n_{\rm b})_{\rm tan} =     \bm \sigma \bm n_{\rm b} - \left( \bm \sigma \bm n_{\rm b} \cdot \bm n_{\rm b} \right) \bm n_{\rm b} &= T_{\rm b} \bB{t}_{\rm b},\label{e:dyn_BC_bottom}\\
(\bm \sigma \bm n_{\rm s})_{\rm tan} =     \bm \sigma \bm n_s - \left( \bm \sigma \bm n_s \cdot \bm n_s \right) \bm n_s &= T_{\rm s} \bB{t}_{\rm s}, \label{e:dyn_BC_surface}
\end{align}
\end{subequations}
where $T_{\rm b}$ and $T_{\rm s}$ are the bottom and surface traction, respectively, and $(\cdot)_{\rm tan}$ denotes the tangent vector. The formulas to compute the bottom traction $T_{\rm b}$ and surface traction $T_{\rm s}$ are given in appendix \ref{sec:stress_tensor_bottom}.

The appropriate boundary conditions in the horizontal coordinate will be discussed later in this article, see Section~\ref{sec:simulations}.

\begin{figure}[t!]

\begin{tabular}{ccc}
(a) Stresses in $z$ space & \hspace{2cm}&
(b) Stresses in $\zeta$ space \\
\includegraphics[scale=0.55]{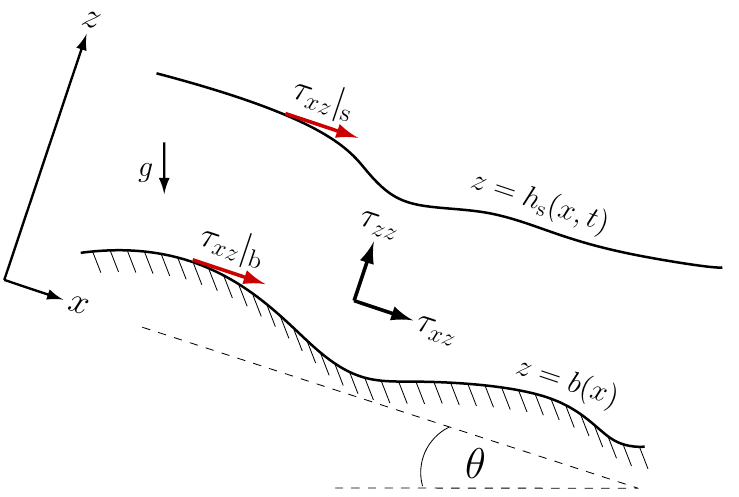} &&
\includegraphics[scale=0.55]{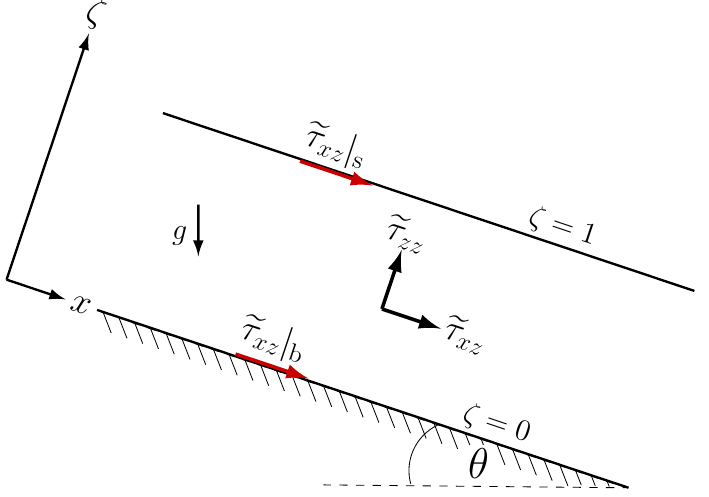}
\end{tabular}
    \caption{Stresses in the physical (left) and transformed (right) coordinate systems.\label{fig:stress}}
\end{figure}

\subsection{Non-dimensional model} \label{sec:non-dim:model}
For a non-dimensional analysis of the model, we introduce a scaling similar to \cite{Garres-Diaz2021}.
As usual for shallow flows, we assume that the ratio of the vertical scale $L$ and the horizontal scale $H$, defined as $\epsilon = H/L$, is a small parameter. Then, defining the factor $U=\sqrt{gL}$ as the characteristic velocity scale, the coordinates and variables in \eqref{eq:govdim} are non-dimensionalized as follows:
\begin{equation} \label{eq:nondim}
\begin{aligned}
  (x,z,t) &= \big(Lx^*, Hz^*, \tfrac{L}{U} t^*\big), \\[1ex]
  (u,w) &= U \left( u^*, \epsilon w^* \right), \\[1ex]
  (p, \tau_{xx},\tau_{zz},\tau_{xz}) &= \rho g\cos(\theta)H\big(p^*, \epsilon \tau_{xx}^*, \epsilon \tau_{zz}^*, \tau_{xz}^*\big), \\[1ex]
  (b,h) &= H(b^*, h^*),
\end{aligned}
\end{equation}
where the asterisk is used to denote the dimensionless variables which are of order $\mathcal{O}(1)$, with respect to the ratio $\epsilon$.
Thus, $\partial_x b = \epsilon \partial_{x^*} b^*$ and $ \partial_x h = \epsilon \partial_{x^*} h^*$ are small variables of order $\epsilon$.
Meanwhile, we see that the symmetric strain-rate tensor scales as
\[
\bB{D}(\bm u) = \frac{U}{H}
\begin{pmatrix}
2\epsilon \partial_{x^*}u^* & \partial_{z^*}u^* + \epsilon^2 \partial_{x^*} w^*\\[1ex]
\partial_{z^*}  u^* + \epsilon^2 \partial_{x^*} w^* & 2\epsilon \partial_{z^*} w^*
\end{pmatrix},
\]
and
\begin{equation} \label{eq:dunorm}
   \| \bB{D}(\bm u) \| = \Big(\tfrac{1}{2}{\rm tr}\big(\bB{D}(\bB{u})^{\rm T}\bB{D}(\bB{u})\big)\Big)^{1/2} = \frac{U}{H} \left | \partial_{z^*} u^*\right | + \mathcal{O}(\epsilon),
\end{equation}
where ${\rm tr}(\cdot)$ represents the trace operator.

The governing incompressible Navier-Stokes equations \eqref{eq:govdim} can then be cast in the following dimensionless form:
\begin{subequations} \label{syst:ndgov}
\begin{align}
  \diffp{u}{x} + \diffp{w}{z} &= 0, \label{eq:ndgova} \\
  \diffp{u}{t} + u \diffp{u}{x} + w \diffp{u}{z} &= \sin(\theta) + \cos(\theta) \diffp{\tau_{xz}}{z} - \epsilon \cos(\theta) \left( \diffp{p}{x} - \epsilon \diffp{\tau_{xx}}{x} \right), \label{eq:ndgovb} \\
  \epsilon \left( \diffp{w}{t} + u \diffp{w}{x} + w \diffp{w}{z} \right) &= -\cos(\theta) -\cos(\theta) \left(\diffp{p}{z} - \epsilon \diffp{\tau_{zz}}{z} - \epsilon \diffp{\tau_{xz}}{x}\right), \label{eq:ndgovc}
\end{align}
\end{subequations}
where the asterisk symbols have been omitted for clarity.
The entries of the deviatoric viscous stress tensor $\bm{\tau}$ are now related to the dimensionless strain rates as follows:
\begin{align}\label{eq:tau:epsilon2}
\begin{pmatrix}
\tau_{xx} & \tau_{xz} \\
\tau_{zx} & \tau_{zz}
\end{pmatrix} = \frac{\eta U}{\rho g \cos(\theta) H^2 }
\begin{pmatrix}
2\partial_x u & \partial_z u + \epsilon^2 \partial_x w \\[1ex]
\partial_z u + \epsilon^2 \partial_x w & 2 \partial_z w
\end{pmatrix}.
\end{align}
Dimensionless system \eqref{syst:ndgov} is accompanied by the respective dimensionless boundary conditions of the original system. Next to the kinematic boundary conditions given by \eqref{eq:kinematic:b} and \eqref{eq:kinematic:hs}, the dimensionless dynamic boundary conditions \eqref{e:dyn_BC_bottom} and \eqref{e:dyn_BC_surface} up to order $\epsilon^2$ are, respectively, given by
\begin{align} \label{e:dyn_BCb_nd}
   (\bm \sigma \bm n_{\rm b})_{\rm tan} & = \tau_{xz}|_b \,\bB{t}_{\rm b},\qquad (\bm \sigma \bm n_{\rm s})_{\rm tan} = - \tau_{xz}|_s \,\bB{t}_{\rm s},
\end{align}
where the bottom and surface stress are denoted by $\tau_{xz}|_b$ and $\tau_{xz}|_s$, respectively.

The long-wave shallowness limit is attained for $\epsilon \to 0$. In this limit $\epsilon \to 0$, equation \eqref{eq:ndgovc} reduces to
\begin{equation} \label{eq:p:hydrostatic}
  p(x,z,t) = b(x) + h(x,t) - z,
\end{equation}
which corresponds to a hydrostatic pressure. All the models studied in this work are based on the assumption that the pressure is hydrostatic. Note that there also exist shallow water models and shallow water moment models for non-hydrostatic pressure \cite{Bassi2020,Garres-Diaz2021,Scholz2023}. We leave further extensions of our approach to non-hydrostatic pressure for future work.

\subsection{Transformed reference system} \label{sec:transformed:ref:syst}
To include information about the vertical variation of the velocity profile, we derive a moment model following the approach in \cite{kowalski2018moment}. For this purpose, we first define a mapping of the vertical coordinate $z \in [b(x),h(x,t)]$ onto a variable with fixed bounds $\zeta \in [0,1]$ via the transformation $\zeta=(z-b(x))/h(x,t)$. Then, we apply this change of variables on each term of system \eqref{syst:ndgov}. In turn, given a scalar function $\psi = \psi(x,z,t)$, the following identities are used to compute the partial derivatives in the new vertical coordinate \cite{kowalski2018moment}:
\begin{align}
 h \diffp{\psi}{z} = \diffp{\widetilde{\psi}}{\zeta}\quad\text{and}\quad h \diffp{\psi}{s} &= \diffp{}{s} \big( h \widetilde{\psi} \big) - \diffp{}{\zeta} \left( \widetilde{\psi} \dfrac{\partial}{\partial s}(\zeta h+b) \right), \quad \text{for } s=x,t,
\end{align}
where $\widetilde{\psi}(x,\zeta,t)=\psi \big(x,h(x,t)\zeta+b(x), t \big)$ is the transformed function. For the transformed hydrostatic pressure, in particular, we have
\begin{equation} \label{eq:pnd}
  \widetilde{p}(x,\zeta,t) = h(x,t)\left( 1 - \zeta \right).
\end{equation}
We can then proceed to transform the remaining equations \eqref{eq:ndgova} and \eqref{eq:ndgovb}.
Following the detailed steps in \cite{kowalski2018moment}, the system of transformed mass conservation and momentum equation is given by 
\begin{subequations}
\begin{align}
  \diffp{}{x} \left( h \widetilde{u} \right) +
\diffp{}{\zeta} \Big( \widetilde{w} - \widetilde{u} \dfrac{\partial}{\partial x} (\zeta h+b) \Big) &= 0, \label{e:ref_mass}\\[1ex]
  \diffp{}{t} (h\widetilde{u}) + \diffp{}{x} \left( h \widetilde{u}^2 + \tfrac{1}{2} \epsilon \cos(\theta) h^2 \right) + \diffp{}{\zeta} \left( \widetilde{u} W  \right)
  &= \sin(\theta) \, h + \cos(\theta)\Big( \diffp{\widetilde{\tau}_{xz}}{\zeta} - \epsilon h \diffp{b}{x}\Big), \label{e:ref_momentum}
\end{align}
\end{subequations}
Note that in \eqref{e:ref_momentum}, the term $\epsilon^2 \partial_x \tau_{xx}$ of order $\epsilon^2$ has been omitted, and the only remaining friction term corresponds to $\widetilde{\tau}_{xz}$ acting in the vertical direction, which is illustrated in Figure~\ref{fig:stress}. Moreover, $W$ in \eqref{e:ref_momentum} corresponds to the so-called vertical coupling term, representing the effect of the vertical velocity $w$ on the $x$-momentum equation \cite{kowalski2018moment}, which is given by:
\begin{align*}
W = \widetilde{w} - \dfrac{\partial}{\partial t}(\zeta h+b) - \widetilde{u} \dfrac{\partial}{\partial x}(\zeta h+b)
= -\diffp{}{x} \left( h \int_0^{\zeta} \widetilde{u} \,{\rm d} \hat{\zeta} \right) + \zeta \dfrac{\partial}{\partial x}(hu_m),
\end{align*}
where the second equality follows from the integrated continuity equation
\begin{align*}
\widetilde{w} = -\diffp{}{x} \left( h \int_0^{\zeta} \widetilde{u}\, {\rm d} \hat{\zeta} \right) +  \widetilde{u} \dfrac{\partial}{\partial x}(\zeta h+b),
\end{align*}
and the mean or averaged velocity $u_m$ is defined as
\begin{align} \label{eq:def:um}
u_m(x,t) = \int_0^1 \widetilde{u}(x,\zeta,t)\,{\rm d}\zeta.
\end{align}

\subsection{Shallow water moment equations} \label{sec:SWME:general}
To obtain the shallow water moment model, the main assumption is that the velocity $\widetilde{u}$ can be written as a polynomial function in the vertical direction, with coefficient functions that depend on the horizontal coordinate $x$ and time $t$ \cite{kowalski2018moment}. In turn, we let $\phi_j=\phi_j(\zeta)$ be the shifted Legendre polynomial of degree $j$ defined for $\zeta\in [0,1]$ and assume that $\widetilde{u}$ is given by a linear combination of $N$ Legendre polynomials plus the depth average $u_m$ as follows:
\begin{equation} \label{eq:mom_ansatz}
  \widetilde{u}(x,\zeta,t)= u_m (x,t) + \sum_{j=1}^N \alpha_j(x,t) \phi_j ( \zeta ),
\end{equation}
where the scalars $\alpha_j$ are the expansion coefficients, also called moments, which are considered unknowns in the model. We remark that $\phi_j(\zeta)=({\rm d}^j/{\rm d}\zeta^j)(\zeta-\zeta^2)^j/j!$ for all $j\in\mathbb{N}$, with $\phi_j(0)=1$ and $\phi_j(1)=(-1)^j$. The integer $N$ is the maximal degree of the expansion, also called the order of the model.

To obtain the conservation of mass equation,  we integrate \eqref{e:ref_mass} over the transformed vertical coordinate, from $\zeta=0$ to $\zeta=1$. Then, using the kinematic boundary conditions \eqref{eq:kinematic_BQ}, we obtain:
\begin{equation} \label{eq:mom_h}
  \diffp{h}{t} + \diffp{}{x}(h u_m) = 0,
\end{equation}
where the mean velocity $u_m$ is defined in \eqref{eq:def:um}. Following the derivation in \cite{kowalski2018moment}, the equation for $u_m$ is obtained from the momentum equation \eqref{e:ref_momentum} after replacing the polynomial expansion \eqref{eq:mom_ansatz} and  integrating over the transformed vertical coordinate $\zeta$, that is
\begin{equation} \label{eq:mom_hu}
  \diffp{}{t} (h u_m) + \diffp{}{x} \bigg( \frac{1}{2} \epsilon \cos(\theta) h^2 + h u_m^2
  + \sum_{j=1}^N \frac{h \alpha_j^2}{2j+1} \bigg) = \sin(\theta) h + \cos(\theta) \bigg( \, \widetilde{\tau}_{xz}\big|_{\rm b}^{\rm s} - \epsilon \, h \diffp{b}{x}\bigg),
\end{equation}
where $\widetilde{\tau}_{xz}\big|_{\rm b}^{\rm s} = \widetilde{\tau}_{xz}|_{\rm s} - \widetilde{\tau}_{xz}|_{\rm b}$ is the difference between the surface friction $\widetilde{\tau}_{xz}|_{\rm s}$ and the bottom friction $\widetilde{\tau}_{xz}|_{\rm b}$. These friction terms will be treated in Section~\ref{sec:framework}.

The equations for the moments $\alpha_i$ are obtained by Galerkin projections of the transformed momentum equation  \eqref{e:ref_momentum} onto Legendre test functions $\phi_i$, such that we obtain for $i=1,\dots,N$:
\begin{equation} \label{eq:mom_hui}
\begin{aligned}
  \diffp{}{t} (h\alpha_i) &+ \diffp{}{x} \bigg( 2hu_m \alpha_i + h \sum_{j,k=1}^N A_{ijk} \alpha_j \alpha_k \bigg) - u_m \diffp{}{x} (h\alpha_i) + \sum_{j,k=1}^N B_{ijk} \alpha_k \diffp{}{x}(h\alpha_j) \\
  &= (2i + 1)\cos (\theta) \, \left( (-1)^i \widetilde{\tau}_{xz}|_{\rm s} - \widetilde{\tau}_{xz}|_{\rm b} -\int_0^1 \phi_i' \widetilde{\tau}_{xz} {\rm d}\zeta \right),
\end{aligned}
\end{equation}
where the coefficients $A_{ijk}, B_{ijk}$ can be precomputed according to the formulas given in \cite{kowalski2018moment} as
\begin{align}
    A_{ijk} &= (2i+1)\int_0^1 \phi_i(\zeta) \phi_j(\zeta) \phi_k(\zeta) \, {\rm d}\zeta,\\
    B_{ijk} &= (2i+1)\int_0^1 \phi_i'(\zeta)\left( \int_0^{\zeta} \phi_j(\hat{\zeta}) \, {\rm d}\hat{\zeta} \right)  \phi_k(\zeta) \, {\rm d}\zeta.
\end{align}
The right hand sides of \eqref{eq:mom_hu} and \eqref{eq:mom_hui} include the basal friction $\widetilde{\tau}_{xz}|_{\rm b}$ and surface friction $\widetilde{\tau}_{xz}|_{\rm s}$, and equation \eqref{eq:mom_hui} also includes the additional bulk friction term $\int_0^1 \phi_i' \widetilde{\tau}_{xz} \,{\rm d}\zeta$. Note that the balance equation for the mean velocity $u_m$ \eqref{eq:mom_hu} does not contain a bulk friction term since the mean velocity is constant throughout the water depth and does not induce bulk shear.
Based on the modelled flow conditions, both basal and bulk friction can be modelled according to different laws as outlined in Section~\ref{sec:framework}.

The moment model with the vector of unknowns $(h,u_m,\alpha_1,\alpha_2,\dots,\alpha_N) \in \mathbb{R}^{N+2}$ is then defined by the coupled system of equations \eqref{eq:mom_h}, \eqref{eq:mom_hu}, and \eqref{eq:mom_hui}. The main difference of this model compared to \cite{kowalski2018moment} is due to the scaling followed to obtain the model equations, which allow us to keep track of the friction term $\tau_{xz}$ and to extend the SWME to general friction models. Furthermore, the model in this paper is derived over spatial coordinates with respect to an inclined plane with angle $\theta\in[0,\tfrac{\pi}{2})$ with respect to the horizontal line, which generalises the classical SWME \cite{kowalski2018moment}.
In the next section, we introduce an alternative moment model which fulfills the hyperbolicity property, that has been developed in \cite{Koellermeier2020}.

\subsection{Hyperbolic shallow water moment equations}

We now consider the vector of conservative variables $\mathbf{U} = (h, h u_m, h \alpha_1, \dots, h \alpha_N)^{\rm T}\in \mathbb{R}^{N+2}$ and the system composed by the conservation of mass \eqref{eq:mom_h}, momentum balance law \eqref{eq:mom_hu}, and moment equations \eqref{eq:mom_hui}. We can also employ the short notation $\mathbf{U} = h(1,u_m,\bB{\alpha})^{\rm T}$, where $\bB{\alpha} = (\alpha_1,\alpha_2,\dots,\alpha_N)^{\rm T}\in \mathbb{R}^N$.
As described in \cite{Koellermeier2020}, the coupled system \eqref{eq:mom_h}, \eqref{eq:mom_hu} and \eqref{eq:mom_hui}, is only conditionally hyperbolic with possible loss of hyperbolicity for certain values of the unknowns.
This limitation lead to instability problems during simulations as shown in particular test cases in \cite{Koellermeier2020}, similar to the loss of hyperbolicity for moment models in rarefied gases, see \cite{Cai2013, Koellermeier2017b, Koellermeier2014, Koellermeier2017}. A solution to overcome this problem is to consider the hyperbolic regularisation of the transport part introduced in \cite{Koellermeier2020} called hyperbolic shallow water moment equations (HSWME).
The hyperbolic model results from a linearisation of the system matrix around linear velocity profiles, i.e. only the effects of the first order moment $\alpha_1$ are fully kept. The model has been studied in detail in, for instance, \cite{Garres2020, KoellermeierQian2022, Pimentel2020}.

The hyperbolic version of the moment model \eqref{eq:mom_h}, \eqref{eq:mom_hu} and \eqref{eq:mom_hui}, in vector notation, is the corresponding HSWME given by
\begin{equation}\label{e:HSWME}
 \dfrac{\partial \mathbf{U}}{\partial t} + \mathbf{A}_{\rm H}(\mathbf{U}) \dfrac{\partial \mathbf{U}}{\partial x} = \mathbf{S}(\mathbf{U}),
\end{equation}
where the system matrix $\mathbf{A}_{\rm H} \in \mathbb{R}^{(N+2)\times(N+2)}$ is defined by
\begin{equation}\label{HSWME_A}
        \mathbf{A}_{\rm H}(\mathbf{U})=\left(
        \begin{array}{ccccccc}
          & 1 &  &  &  &  &  \\
        \epsilon \cos(\theta) h-u_m^2-\frac{1}{3} \alpha_1^2& 2 u_m & \frac{2}{3} \alpha_1 &  &  &  &  \\
         -2u_m  \alpha_1 & 2 \alpha_1 & u_m & \frac{3}{5} \alpha_1  & &  \\
         -\frac{2}{3} \alpha_1^2 &  & \frac{1}{3} \alpha_1 & u_m & \ddots  &  \\
         &  &  & \ddots & \ddots  & \frac{N+1}{2N+1} \alpha_1  \\
          &  &  &  & \frac{N-1}{2N-1} \alpha_1 & u_m  \\
        \end{array}
        \right),
\end{equation}
with all other entries equal to zero. The right-hand side source term vector $\mathbf{S}(\mathbf{U}) \in \mathbb{R}^{N+2}$ is not impacted by the hyperbolic regularisation as it
models friction and external effects such as the bottom topography and is given by
\begin{equation}\label{e:friction_terms}
\mathbf{S}(\mathbf{U}) = \left(
        \begin{array}{c}
            0 \\[1ex]
            \sin (\theta) h + \cos(\theta) \big( \widetilde{\tau}_{xz}\big|_{\rm b}^{\rm s} - \epsilon h \partial_x b \big)\\[1ex]
            (2\cdot 1+1) \cos(\theta) \left( (-1)^1\widetilde{\tau}_{xz}|_{\rm s} - \widetilde{\tau}_{xz}|_{\rm b} - \int_0^1 \phi_1' \widetilde{\tau}_{xz} \, {\rm d}\zeta \right) \\[1ex]
            \vdots \\[1ex]
            (2\cdot N + 1)\cos (\theta)  \left( (-1)^N\widetilde{\tau}_{xz}|_{\rm s} - \widetilde{\tau}_{xz}|_{\rm b} - \int_0^1 \phi_N' \widetilde{\tau}_{xz}\, {\rm d}\zeta \right)
        \end{array}
        \right),
\end{equation}
where the friction term $\widetilde{\tau}_{xz}$ is still to be determined, specifically its definition at the bottom $\zeta = 0$, at the top $\zeta = 1$, and within the values $\zeta\in[0,1]$ for the integral terms in \eqref{e:friction_terms}. Both models, the general moment model given by equations \eqref{eq:mom_h}, \eqref{eq:mom_hu} and \eqref{eq:mom_hui}, and the HSWME model \eqref{e:HSWME}, are part of a hierarchical modelling framework to derive model hierarchies for general shallow flows as an extension of the shallow water moment model.

\begin{remark}
In general, one-dimensional shallow water moment models can be written as the following first order system:
\begin{equation}\label{FBSsystem}
\dfrac{\partial \mathbf{U}}{\partial t} + \dfrac{\partial \mathbf{F}(\mathbf{U})}{\partial x} + \mathbf{A}(\mathbf{U}) \dfrac{\partial \mathbf{U}}{\partial x} = \mathbf{S}(\mathbf{U}) := \mathbf{S}_0(\mathbf{U})+\mathbf{S}_1(\mathbf{U}) \dfrac{\partial b}{\partial x},
\end{equation}
where $\mathbf{F}\in\mathbb{R}^{N+2}$ is called the conservative flux, $\mathbf{A}\in \mathbb{R}^{(N+2)\times (N+2)}$ is the matrix of non-conservative fluxes, $\mathbf{S}_1\in \mathbb{R}^{N+2}$ is a vector containing the factors of $\partial_x b$ in the source term, and $\mathbf{S}_0\in \mathbb{R}^{N+2}$ is a source term (nonlinear friction term). The given form of system \eqref{e:HSWME}, can be obtained using $\mathbf{F} = \mathbf{0}$, and $\mathbf{A} = \mathbf{A}_{\rm H}$ and
\begin{align*}
 \mathbf{S}_0(\mathbf{U}) & = \cos(\theta) \Big(\mathbf{S}_{\rm lim}(\mathbf{U}) - \mathbf{S}_{\rm bulk}(\mathbf{U})\Big) + \sin(\theta)\mathbf{S}_{\rm g}(\mathbf{U}),\\[1ex]
 \mathbf{S}_1(\mathbf{U}) & = -\cos(\theta)\big(0,\epsilon h, 0,\dots,0\big)^{\rm T},
\end{align*}
where $\mathbf{S}_{\rm lim}$ is the source term related to the friction at the top and bottom surfaces, $\mathbf{S}_{\rm bulk}$ contains the bulk friction terms, and $\mathbf{S}_{\rm g}$ is the gravity related source, which only affects the mean average velocity equation. These friction terms are defined componentwise by
\begin{align*}
\begin{array}{lll}
 \mathrm{S}_{{\rm lim},1} = 0, & \qquad \mathrm{S}_{{\rm lim},2} = \widetilde{\tau}_{xz}|_{\rm b}^{\rm s},
 &\qquad \mathrm{S}_{{\rm lim},i+2} = (2i+1)\big((-1)^{i} \widetilde{\tau}_{xz}|_{\rm s} - \widetilde{\tau}_{xz}|_{\rm b}\big),\\[1ex]
  \mathrm{S}_{{\rm bulk},1} = 0,& \qquad \mathrm{S}_{{\rm bulk},2}  = 0,  &\qquad  \mathrm{S}_{{\rm bulk},i+2} = (2i+1)\mathcal{T}_{i},\\[1ex]
 \mathrm{S}_{{\rm g},1} = 0,&\qquad  \mathrm{S}_{{\rm g},2}  = h,  &\qquad \mathrm{S}_{{\rm g},i+2}  = 0,
\end{array}
\end{align*}
for $i=1,\dots, N$, and the additional bulk friction term $\mathcal{T}_i$ is defined as
\begin{align}\label{eq:taubulk:integral}
\mathcal{T}_{i}(x,t) := \int_0^1 \widetilde{\tau}_{xz}(x,\zeta,t) \phi_i'(\zeta)\, {\rm d}\zeta.
\end{align}

\end{remark}

From now on, and without loss of generality, we are going to consider the HSWME model instead of the general SWME from Section~\ref{sec:SWME:general}. Furthermore, all numerical examples shown in Section~\ref{sec:simulations} are computed using the HSWME given by \eqref{e:HSWME}. The focus is therefore on the right-hand side friction terms.

\section{Friction models} \label{sec:framework}

In this section, we derive a 3-step framework that includes general friction terms depending on the stress tensor entry $\tau_{xz}$, which can be modelled by different expressions depending on the type of rheology. For the moment models presented in Section \ref{sec:SWmodels}, we observe that $\tau_{xz}$ only appears on the right-hand side of the equations, while the transport part of the models remain independent of the friction.
Therefore, we turn our attention primarily to the right-hand side friction terms and, without loss of generality, will consider the HSWME model \eqref{e:HSWME} from now on.

For the source term in \eqref{e:friction_terms}, the transformed friction $\widetilde{\tau}_{xz}$ is evaluated at the bottom $\zeta=0$ and free surface $\zeta=1$, and is also present in the additional bulk friction $\mathcal{T}_i$ \eqref{eq:taubulk:integral},
for which we also require a precise definition of $\widetilde{\tau}_{xz}$ for all values $\zeta\in(0,1)$. While a whole friction model is necessary, the friction processes at the bottom and in the interior of the water column can in principle be different. Therefore, we will now present friction models taking this into account. To properly model the friction term, we propose the following 3-step model framework:
\begin{enumerate}
    \item Define the (dimensional) physical constitutive laws for the stress tensor entry $\tau_{xz}$ at the bottom boundary, inside the bulk fluid, and at the surface boundary:
    \begin{align}
        \tau_{xz}(x,z,t) = \begin{cases}
                       \tau_{xz}|_{\textrm{s}}(x,t), & z = h_{\rm s}(x,t), \\[1ex]
                       \tau_{xz}|_{\textrm{bulk}}(x,z,t), & z \in \big(b(x),h_{\rm s}(x,t)\big), \\[1ex]
                       \tau_{xz}|_{\textrm{b}}(x,t), & z = b(x).
                       \end{cases}  \label{eq:tau_general}
    \end{align}
    \item Non-dimensionalise $\tau_{xz}$ following the scaling in \eqref{eq:nondim}, and apply the mapping $z\mapsto \zeta$ described in Section~\ref{sec:transformed:ref:syst}:
    \begin{equation}
        \widetilde \tau_{xz}(x,\zeta,t) = \frac{1}{\rho g \cos (\theta) H} \tau_{xz}\big(x,h\zeta+b,t\big),\label{eq:general:tildetau}
    \end{equation}
where $h$ and $b$ are understood to be dimensionless variables.
    
    \item Evaluation of the dimensionless and transformed friction terms in \eqref{e:friction_terms} by computing the respective boundary terms at $\zeta=0$ and $\zeta=1$, and integral term \eqref{eq:taubulk:integral} to derive the hierarchical model:
    \begin{equation}\label{e:boundary_interior_terms}   
        \widetilde{\tau}_{xz}|_{\rm b}, \quad
        \widetilde{\tau}_{xz}|_{\rm s}, \quad \text{and} \quad
       \mathcal{T}_i = \int_0^1 \widetilde{\tau}_{xz}(x,\zeta,t) \phi_i'(\zeta)\, {\rm d}\zeta, \quad \text{for }i=1,\ldots,N.
    \end{equation}
\end{enumerate}
The three consecutive steps constitute a clearly defined procedure to derive the source term for any friction model.
The concise and general model framework has three major benefits. Firstly, due to the concise definition and the usage of the general shallow flow model, the framework is widely applicable. Secondly, most of the derived models are given analytically, which allows for analysis of the models without further numerical approximation or simulation. This distinguishes our framework from existing multilayer or heuristic models like \cite{FernandezNieto2016}. Thirdly, the splitting of the friction term into boundary and internal friction makes it possible to separate the two physically different effects, which is important to attach physical meaning and interpretation to the resulting model.

\begin{remark}
    Most of the following models are directly given in analytical form. While the evaluation of the boundary terms at the bottom and at the free surface are straightforward, the integral term for the interior friction \eqref{e:boundary_interior_terms}, might require a numerical integration procedure, depending on the complexity of the friction model.
\end{remark}

In the upcoming subsections, we apply the described framework \eqref{eq:general:tildetau}--\eqref{e:boundary_interior_terms} to several existing friction models from the literature, while in Section~\ref{sec:granular}, we use this procedure to develop a new model for shallow granular flows based on the $\mu(I)$-rheology.
Table~\ref{table:models} summarizes the friction terms considered in the upcoming section. For clarity, in the remainder of the paper, we omit the dependence of the variables on $x$ and $t$, and focus on the $\zeta$-coordinate whenever possible.

\begin{table}
\renewcommand{\arraystretch}{1.3}
\begin{center}
\begin{tabular}{ c|l|l|l }
 Section &
 \multicolumn{1}{c|}{Model} &
 \multicolumn{1}{c|}{bulk stress $\tau_{xz}(z)$} &
 \multicolumn{1}{c}{bottom stress $\tau_{xz}|_b$} \\[1ex] \specialrule{.1em}{0em}{0.01em}
 Section \ref{sec:Newtonian_Slip_Flow} & Newtonian slip\cite{kowalski2018moment} & Newtonian & slip \\ 
 Section~\ref{sec:Newtonian_Manning} & Newtonian Manning \cite{Garres2020} & Newtonian & Manning-type \\ 
 Section~\ref{sec:Savage-Hutter_model} & Savage-Hutter \cite{Garres-Diaz2021} & Mohr-Coulomb & solid law with local angle $\delta$\\
 Section~\ref{sec:Coulomb-type_const} & Coulomb-type & Mohr-Coulomb & solid law with local angle $\delta$\\
 Section~\ref{sec:granular} & $\mu(I)$-rheology slip & $\mu(I)$-rheology & slip\\
 \bottomrule
\end{tabular}
\end{center}
\caption{\label{table:models}Existing and new models derived using the new framework.}
\end{table}

\subsection{Newtonian slip flow} \label{sec:Newtonian_Slip_Flow}
For the Newtonian bulk flow, which has been used in \cite{kowalski2018moment}, the constitutive law is consistent with \eqref{eq:constit} using a constant viscosity $\eta\geq 0$. In this case, we assume a stress-free surface and constant slip length $\Lambda> 0$ at the bottom.
After neglecting the term $\epsilon^2\partial_x w$ arising in \eqref{eq:tau:epsilon2} to obtain $\tau_{xz}$ in the bulk, and considering slip condition for the bottom friction, the three-step framework reads as follows:
\begin{enumerate}
    \item Dimensional stress in $z$-coordinates:
    \begin{align}
     \tau_{xz} |_{\rm b} = \frac{\eta}{\Lambda} u(b),\qquad \tau_{xz} |_{\rm bulk} = \eta \partial_z u(z),\qquad  \tau_{xz} |_{\rm s} = 0.
     \label{eq:newton}
    \end{align}
    \item Dimensionless stress in $\zeta$-coordinates:
    \begin{align}
     \widetilde{\tau}_{xz}|_{\rm b} = \frac{\nu}{\lambda} \widetilde u(0),\qquad \widetilde{\tau}_{xz}|_{\rm bulk} = \frac{\nu}{h} \partial_{\zeta} \widetilde u(\zeta),\qquad  \widetilde{\tau}_{xz}|_{\rm s} = 0,
\end{align}
where $\lambda = \Lambda /H$ is the dimensionless slip length and $\nu = \eta U / (\rho g \cos(\theta)  H^2)$ is the dimensionless viscosity.
    \item Boundary terms and additional bulk friction:
    \begin{align}
      &  \widetilde{\tau}_{xz}|_{\rm b} = \frac{\nu}{\lambda} \left( u_m + \sum_{i=1}^N \alpha_i \right),\qquad \mathcal{T}_i = \frac{\nu}{h} \sum_{j=1}^N C_{ij} \alpha_j,  \qquad \widetilde{\tau}_{xz}|_{\rm s} = 0,\label{eq:NS_tau_integral}
    \end{align}
    for  $i=1,\ldots,N$ where the coefficients $C_{ij}$ are defined as \cite{kowalski2018moment}
\begin{align}
C_{ij} = \int_{0}^{1} \partial_{\zeta} \phi_i(\zeta) \, \partial_{\zeta} \phi_j(\zeta) \,{\rm d}\zeta,\qquad \text{for } i,j=1,\dots,N. \label{eq:defn:Cij}
\end{align}

\end{enumerate}

\begin{remark}
Note that in this case, the bottom friction $\widetilde{\tau}_{xz}|_{\rm b}$ depends on the bottom velocity $\widetilde{u}(0)$, which thanks to the property of the basis functions $\phi_i(0) = 1$ for all $i=1,\dots,N$, can be written as
\begin{align}\label{eq:bottom:velocity}
\widetilde{u}(0) = u_m + \sum_{i=1}^N \alpha_i.
\end{align}
This equation will be used to compute the bottom velocity in the examples shown in Section~\ref{sec:simulations}.
\end{remark}

\subsection{Newtonian Manning} \label{sec:Newtonian_Manning}
We now consider the friction term proposed in \cite{Garres2020} which was applied to a shallow water moment model for bedload sediment transport. The model was derived by using a Newtonian bulk flow and a Manning bottom stress term.  In comparison with the Newtonian slip flow model in Section~\ref{sec:Newtonian_Slip_Flow}, the main change is the bottom friction term $\tau_{xz}|_{\rm b}$ which is driven by a Manning law. The stress-free surface and the interior Newtonian friction remain the same. Letting $n\geq0$ be the Manning coefficient, this model can conveniently be written in our framework as follows:
\begin{enumerate}
    \item Dimensional stress in $z$-coordinate:
        \begin{align}
     \tau_{xz} |_{\rm b} = \frac{\rho g n^2}{h^{1/3}} u(b) | u(b) |,\qquad \tau_{xz} |_{\rm bulk} = \eta \partial_z u(z),\qquad  \tau_{xz} |_{\rm s} = 0.\label{eq:new-man}
\end{align}
    \item Dimensionless stress in $\zeta$-coordinate:
\begin{align}
     \widetilde{\tau}_{xz}|_{\rm b} = \frac{\mathfrak{n}^2}{h^{1/3}} \widetilde u(0) | \widetilde u(0) |,\qquad \widetilde{\tau}_{xz}|_{\rm bulk} = \frac{\nu}{h} \partial_{\zeta} \widetilde u(\zeta),\qquad  \widetilde{\tau}_{xz}|_{\rm s} = 0,
\end{align}
  where $\mathfrak{n}^2=n^2 U^2 /(H^{4/3}\cos(\theta))$ and $\nu = \eta U / (\rho g \cos (\theta) H^2)$.
\item Boundary terms and additional bulk friction:
  \begin{align}
    \widetilde{\tau}_{xz}|_{\rm b} =
    \frac{\mathfrak{n}^2}{h^{1/3}} \left| u_m + \sum_{i=1}^N \alpha_i \right| \left( u_m + \sum_{i=1}^N \alpha_i \right),
        \quad \mathcal{T}_i = \frac{\nu}{h} \sum_{j=1}^N C_{ij} \alpha_j,\quad \widetilde{\tau}_{xz}|_{\rm s} = 0,
  \end{align}
for $i = 1,\dots,N$, where $C_{ij}$ is defined in \eqref{eq:defn:Cij}.
\end{enumerate}

\subsection{Savage-Hutter model}
\label{sec:Savage-Hutter_model}
The Savage-Hutter model is a well-known approach to model granular flows like landslides and avalanches \cite{Savage1989,Savage1991}. The friction term described by this model is based on a \textit{rate-independent} assumption, where the shear stress is related to the normal stress via a constant internal friction angle $\varphi$ and a local bed friction angle $\delta$. The bed friction angle satisfies that $\delta <  \theta$, while the internal friction angle is such that $\varphi \geq \delta$. This type of friction model can be included in our concise framework as follows:
\begin{enumerate}
    \item Dimensional stress in $z$-coordinate:
    \begin{align}
    \tau_{xz}|_{\rm b} =  p(b)\sgn(u(b)) \tan(\delta) , \quad \tau_{xz}|_{\rm bulk} = p(z)\sgn(u(z)) \tan(\varphi), \quad \tau_{xz}|_{\rm s} = 0.
    \end{align}
\item Dimensionless stress in $\zeta$-coordinate:
    \begin{align}
    \widetilde{\tau}_{xz}|_{\rm b} = h \sgn(\widetilde u(0)) \tan(\delta), \quad
    \widetilde{\tau}_{xz}|_{\rm bulk} = h(1-\zeta) \sgn(\widetilde u(\zeta)) \tan(\varphi), \quad
    \widetilde{\tau}_{xz}|_{\rm s} = 0,
    \end{align}
where we have used \eqref{eq:pnd} for the pressure $p(z)$. 

\item Boundary terms and additional bulk friction, under the assumption $u(\zeta)>0$ for all $0\leq \zeta \leq 1$:
\begin{align}\label{eq:SHbulk}
\begin{aligned}
\widetilde{\tau}_{xz}|_{\rm b} = h\tan(\delta), \qquad \mathcal{T}_i  = -h\tan(\varphi),\qquad \widetilde{\tau}_{xz}|_{\rm s} = 0,
\end{aligned}
\end{align}
for $i=1,\ldots,N$.
\end{enumerate}
The positive velocity condition $u(\zeta)>0$ is common for avalanches. The calculation of the additional bulk friction term $\mathcal{T}_i$ in
 \eqref{eq:SHbulk} results from an integration by parts and the properties of the Legendre polynomials $\phi_i(\zeta)$ for $i\ge 1$
\begin{equation*}
    \int_0^1 \phi_i'(\zeta) (1-\zeta)\,{\rm d}\zeta = \big(\phi_i(\zeta)(1-\zeta)\big)\Big|_{\zeta=0}
    ^{\zeta= 1} + \int_0^1 \phi_i(\zeta) \,{\rm d}\zeta = -1.
\end{equation*}

Note how the boundary and interior friction terms only differ in sign and the different friction angles. Importantly, different values of the internal friction angle $\varphi$ and the local bed friction angle $\delta$ ensure non-zero total friction that eventually relaxes the higher order moments. In addition, we observe that all moment equations exhibit the same interior friction term given by \eqref{eq:SHbulk}. This is due to the specific formulation with constant friction angle $\varphi$. One way to mitigate this nonphysical behaviour is to use varying friction angles, outlined in the next section.

\subsection{Coulomb-type friction with a constant coefficient}\label{sec:Coulomb-type_const}

For dense granular flows, the constitutive equation for the viscous-stress tensor \eqref{eq:constit} is widely adopted in modern literature \cite{Garres2020}. In this case, the (dimensional) effective viscosity $\eta$ depends on the pressure $p$ and the norm of the symmetric strain-rate tensor $\| \bm D(\bm u)\|$ as
\begin{equation} \label{eq:colvis}
    \eta = \mu\frac{ p}{\|\bm D(\bm u)\|},
\end{equation}
where the constant $\mu$ is a dimensionless friction coefficient. The case of variable $\mu$ is treated in Section \ref{sec:granular} in the context of granular flows. Note that thanks to the hydrostatic-pressure assumption $p$ can be written in terms of $b$, $h$ and $z$ through Equation \eqref{eq:p:hydrostatic}. Furthermore, neglecting the $\mathcal{O}(\epsilon)$ terms in \eqref{eq:dunorm}, we can approximate the dimensional version of $\|\bB{D}(\bB{u})\|$ by $|\partial_z u|$. Thus, replacing the effective viscosity function \eqref{eq:colvis} in the stress-tensor \eqref{eq:constit}, and neglecting the term $\partial_x w$ of order $\mathcal{O}(\epsilon^2)$, we have
\begin{equation}\label{eq:Coulomb:tau}
    \tau_{xz} = \mu p(z) \sgn(\partial_z u), \qquad \text{for }b \leq  z \leq  h_{\rm s}.
\end{equation}
Then, we can straightforwardly obtain the dimensionless transformed stress
\begin{equation} \label{eq:Coulomb:dltau}
    \widetilde{\tau}_{xz} = \mu h(1-\zeta)\sgn(\partial_\zeta \widetilde{u}), \qquad \text{for }0\leq \zeta\leq 1,
\end{equation}
where the factor $\rho g \cos(\theta) H$ arising in the dimensionless form of $\tau_{xz}$ cancels with the one from the dimensionless version of $p$.
We observe that due to the term $\sgn(\partial_\zeta \widetilde u)$ in \eqref{eq:Coulomb:dltau}, an analytical expression for the additional bulk stress $\mathcal{T}_i$ may not exist in general for all $i=1,\dots,N$. Nonetheless, it seems physically reasonable to assume an increasing velocity profile with respect to $\zeta$ such that $\sgn(\partial_{\zeta}\widetilde{u})= 1$ in \eqref{eq:Coulomb:dltau}.

Then, considering a Coulomb law for the bottom friction with constant friction angle $\delta$, the friction model can be described as follows:

\begin{enumerate}
\item Dimensional stress in $z$-coordinate:
\begin{alignat*}{5}
&&
\tau_{xz}|_{\rm b}    &=  p(b) \sgn(u(b)) \tan(\delta),&&\qquad
\tau_{xz}|_{\rm bulk} &\,= \mu p(z) \sgn(\partial_z u), &&\qquad
\tau_{xz}|_{\rm s}    &= 0.
\end{alignat*}

\item Dimensionless stress in $\zeta$-coordinate:
\begin{alignat*}{5}
&&
\widetilde{\tau}_{xz}|_{\rm b}    &=  h\sgn(\widetilde{u}(0))\tan(\delta), &&\qquad
\widetilde{\tau}_{xz}|_{\rm bulk} &\,= \mu h(1-\zeta) \sgn(\partial_\zeta \widetilde{u}), &&\qquad
\widetilde{\tau}_{xz}|_{\rm s}    &= 0.
\end{alignat*}
\item Boundary terms and additional bulk friction, under the assumption $\partial_\zeta \widetilde{u} >0$ for all $0\leq \zeta\leq 1$:
\begin{align}\label{eq:tau:bottom:bulk:Coulomb-type}
\widetilde{\tau}_{xz}|_{\rm b} = h \sgn (\widetilde{u}(0) )\tan(\delta),\qquad \mathcal{T}_i = -\mu h,\qquad \widetilde{\tau}_{xz}|_{\rm s} = 0,
\end{align}
for $i=1,\dots,N$.

\end{enumerate}

\begin{remark}
 Note that the Savage-Hutter friction terms $\widetilde{\tau}_{xz}|_{\rm b}$ and $\mathcal{T}_i$ under the assumption that $\widetilde{u}(\zeta)>0$ are equivalent to the Coulomb-type friction terms when $\partial_\zeta \widetilde{u} >0$ and $\mu = \tan(\varphi)$.
\end{remark}

In the next section, we are going to study more in detail a friction model based on the so-called $\mu(I)$-rheology \cite{Jop2006} typically employed to describe sliding Granular material.

\section{Granular shallow flow equations with $\mu(I)$-rheology} \label{sec:granular}
In this section, we derive models based on the shallow moment equations for granular flows in inclined planes for complex rheologies. This is motivated by the fact that both the Newtonian fluid assumption in Sections \ref{sec:Newtonian_Slip_Flow} and \ref{sec:Newtonian_Manning},  and constant friction coefficient in Section \ref{sec:Coulomb-type_const} do not reflect the complexity of realistic fluids \cite{Barker2017,Forterre2008,Zhuang2025}. We therefore use a variable friction coefficient as suggested in \cite{Jop2006}:
\begin{equation} \label{eq:colvis_muI}
    \eta = \mu(I) \frac{p}{\|\bm D(\bm u)\|}\qquad\text{with}\qquad \mu(I) = \mu_{\rm s} + \frac{\mu_2 - \mu_{\rm s}}{I_0+I}I,
\end{equation}
where $\mu_2,\mu_{\rm s},I_0$ are constants and the inertial number $I$ is given by
\begin{equation} \label{eq:Idef}
    I=\frac{d_{\rm s} \|\bm D(\bm u)\|}{\sqrt{p/\rho_{\rm s}}},
\end{equation}
with $d_{\rm s}$ and $\rho_{\rm s}$ denoting the diameter and density of the granular particles, respectively.
Following the approximations made in Section~\eqref{sec:non-dim:model}, in the same way as in previous sections, we use the approximation $ \|\bm D(\bm u)\| \approx |\partial_z u|$ such that $I = (d_{\rm s}|\partial_z u|) /\sqrt{p/\rho_{\rm s}}$, and then we have
\begin{align*}
    \eta = \mu(I) \frac{p}{|\partial_z u |}, \qquad \mu(I) =
     \mu_{\rm s} + \frac{\mu_2 - \mu_{\rm s}}{  \sqrt{p/\rho_{\rm s}}I_0/ d_{\rm s} + |\partial_z u|}|\partial_z u|.
\end{align*}
Then, following the bulk definition $\tau_{xz} = \eta\partial_z u$, we can determine the $\mu(I)$-rheology based friction terms. For the bottom friction, it is possible to consider a slip condition as in Section~\ref{sec:Newtonian_Slip_Flow}, a Manning law as in Section~\ref{sec:Newtonian_Manning} or a Coulomb based friction as in Section~\ref{sec:Savage-Hutter_model}. For the purpose of presentation, we first consider a slip condition for the bottom friction and refer to other bottom friction laws in Remark~\ref{r:muIbottomfriction}. The friction terms in this case are:
\begin{enumerate}
    \item Dimensional stress in $z$-coordinate:
    \begin{align}
\tau_{xz}|_{\rm b} = \frac{\eta_0}{\Lambda} u(b),\quad
\tau_{xz}|_{\rm bulk} = \bigg(\mu_{\rm s} + \frac{(\mu_2 - \mu_{\rm s})|\partial_z u|}{  c_I\sqrt{p} + |\partial_z u|}\bigg)p\,\sgn(\partial_z u), \quad \tau_{xz}|_{\rm s} = 0,\label{eq:tauxz_muI}
    \end{align}
where $c_I = I_0/ (d_{\rm s}\sqrt{\rho_{\rm s}})$, and $\eta_0$ is a viscosity parameter.
    \item Dimensionless stress in $\zeta$-coordinate:
        \begin{align}
        \widetilde{\tau}_{xz}|_{\rm b} = \frac{\nu_0}{\lambda} \widetilde u(0),\quad
       \widetilde{\tau}_{xz}|_{\rm bulk} = \bigg(\mu_{\rm s} + \frac{(\mu_2 - \mu_{\rm s})|\partial_\zeta \widetilde{u}|}{ \tilde{c}_I h^{3/2}\sqrt{1-\zeta} + |\partial_\zeta \widetilde{u}|}\bigg)h(1-\zeta)\,\sgn(\partial_\zeta \widetilde{u}), \label{eq:tauxz_muI2}
    \end{align}
where $\tilde{c}_I = (I_0H/d_{\rm s})\big((\rho/\rho_{\rm s}) \epsilon \cos(\theta)\big)^{1/2}$, $\lambda = \Lambda /H$, $\nu_0 = \eta_0 U / (\rho g \cos(\theta)  H^2)$ and $\widetilde{\tau}_{xz}|_{\rm s} = 0$.
 \item Boundary terms and additional bulk friction:
    \begin{align}
        & \widetilde{\tau}_{xz}|_{\rm s} = 0, \qquad
       \widetilde{\tau}_{xz}|_{\rm b} = \frac{\nu_0}{\lambda} \left( u_m + \sum_{i=1}^N \alpha_i \right), \label{eq:taub_allgr}\\[1ex]
       & \mathcal{T}_i = h \int_0^1 \phi_i'(\zeta) (1-\zeta) \bigg(\mu_{\rm s} + \frac{(\mu_2 - \mu_{\rm s})|\partial_\zeta \widetilde{u}|}{ \tilde{c}_I h^{3/2}\sqrt{1-\zeta} + |\partial_\zeta \widetilde{u}|}\bigg) \sgn(\partial_{\zeta}\widetilde{u}) \,{\rm d}\zeta, \label{eq:bulkmuI}
    \end{align}
for $i=1,\dots,N$.
\end{enumerate}

In general, the main difficulty in obtaining a closed expression for the integral in \eqref{eq:bulkmuI} relies in the fact that $\partial_\zeta \widetilde{u}(\zeta)$, which is a polynomial function in $\zeta$, may vary in its sign.
However, we can split the interval $(0,1)$ into the union of two subintervals with positive and negative derivative of $\widetilde{u}$, respectively, therefore we let $\mathcal{J}^+,\mathcal{J}^{-}\subseteq (0,1)$ be the subintervals defined as:
\begin{align*}
 \zeta\in \mathcal{J}^{+}: \quad \partial_\zeta \widetilde{u}(\zeta)>0\qquad \text{and}\qquad  \zeta\in \mathcal{J}^{-}: \quad \partial_\zeta \widetilde{u}(\zeta)<0.
\end{align*}
Thus, the bulk friction term in \eqref{eq:bulkmuI} can be rewritten as
\begin{align*}
  \mathcal{T}_i & = h  \int_{\mathcal{J}^+}\phi_i'(\zeta)  (1-\zeta) \left( \mu_s + \frac{(\mu_2-\mu_s) \partial_\zeta \widetilde u(\zeta)}{\partial_\zeta \widetilde u(\zeta) + \tilde{c}_I h^{3/2} \sqrt{1-\zeta}} \right) {\rm d}\zeta\\
  &\quad -  h  \int_{\mathcal{J}^-}\phi_i'(\zeta)  (1-\zeta) \left( \mu_s + \frac{(\mu_2-\mu_s) \partial_\zeta \widetilde u(\zeta)}{\partial_\zeta \widetilde u(\zeta) - \tilde{c}_I h^{3/2} \sqrt{1-\zeta}} \right) {\rm d}\zeta, \quad\text{for } i = 1,\dots,N,
\end{align*}
Now, using the change of variables $\zeta = 1- \xi^2 $ such that ${\rm d}\zeta = -2 \xi {\rm d}\xi$ in both integrals above, we get the following expression:
\begin{equation} \label{eq:bulkmuI_N}
    \mathcal{T}_i = 2h \bigg( \int_{\mathcal{J}^+} \mathcal{G}^+_i (\xi) \,{\rm d}\xi - \int_{\mathcal{J}^{-}} \mathcal{G}^-_i(\xi)\, {\rm d}\xi\bigg),
\end{equation}
where
\begin{equation*}
    \mathcal{G}^\pm_i(\xi) = \xi^3 q_i(\xi)\left( \mu_s + \frac{(\mu_2-\mu_s) \sum_{j=1}^N \alpha_jq_j(\xi)}{\sum_{j=1}^N \alpha_jq_j(\xi) \,\pm\, \tilde{c}_I h^{3/2}\xi} \right),\quad\text{with}\quad
    q_i(\xi) = \phi'_i(1-\xi^2).
\end{equation*}
Note that for all $i=1,\dots,N$, $q_i$ is a polynomial function, and then both $\mathcal{G}_i^{+}$ and  $\mathcal{G}_i^{-}$ are rational functions so that the integral defining $\mathcal{T}_i$ can, in principle, be computed analytically. However, in the numerical examples shown in Section~\ref{sec:simulations}, we approximate $\mathcal{T}_i$ through a high order Gauss-Legendre quadrature rule for simplicity. 

\begin{remark}\label{r:muIbottomfriction}
In this section we have chosen a slip condition for the bottom friction. However, the choice of the bottom friction $\widetilde{\tau}_{xz}|_{\rm b}$ may not be restricted to only slip condition. Some of the possible choices for the bottom friction are: 
\begin{align}
\text{Manning-law:} \hspace{-1.5cm}&& \widetilde{\tau}_{xz}|_{\rm b} & =\frac{\mathfrak{n}^2}{h^{1/3}} \widetilde{u}(0)\left| \widetilde{u}(0) \right|, \label{eq:}\\
\text{Coulomb-law:}\hspace{-1.5cm}&& \widetilde{\tau}_{xz}|_{\rm b} & =h \sgn (\widetilde{u}(0) )\tan(\delta),\\
\text{$\mu(I)$-rheology:}\hspace{-1.5cm}&& \widetilde{\tau}_{xz}|_{\rm b} & 
= \bigg(\mu_{\rm s} + \frac{(\mu_2 - \mu_{\rm s})|\partial_\zeta \widetilde{u}(0)|}{ \tilde{c}_I h^{3/2} + |\partial_\zeta \widetilde{u}(0)|}\bigg)h\,\sgn(\partial_\zeta \widetilde{u}(0)), \label{eq:muI:bottom}
\end{align}
where $\delta$ and $\mathfrak{n}$ are explained in Section~\ref{sec:Newtonian_Manning} and \ref{sec:Savage-Hutter_model}, respectively. The last friction term \eqref{eq:muI:bottom} corresponds to the case of assuming also a $\mu(I)$-rheology friction at the bottom, which for flows starting from rest may lead to zero variations in the moments. In Section~\ref{sec:equilibrium:stability}, we provide an analysis of equilibrium states for the case $N=1$ of the $\mu(I)$-rheology bottom friction considering $\widetilde{\tau}_{xz}|_{\rm b}$ as in \eqref{eq:muI:bottom}.
\end{remark}
In the next section, we exemplify the model for $N=1$, where a linear approximation for the velocity is assumed, and provide an analytical expression for the bulk friction term \eqref{eq:bulkmuI}.

\subsection{Case $N=1$}

We now focus exemplarily on the case when $N=1$, i.e.~for linear variations of the velocity $\widetilde{u}$ in the transformed vertical coordinate $\zeta$, see \eqref{eq:mom_ansatz}. In this case, the only moment is $\alpha_1$ and the vector of unknowns becomes $\mathbf{U} = (h,hu_m,h\alpha_1)^{\rm T}$. Furthermore, the system matrix $\mathbf{A}_{\rm H}$ and source term $\mathbf{S}$ in the HSWME model \eqref{e:HSWME} are reduced to
\begin{align*}
\mathbf{A}_{\rm H}(\mathbf{U}) =
\begin{pmatrix}
0 & 1 & 0 \\
\epsilon \cos(\theta) h-u_m^2-\frac{1}{3}\alpha_1^2 & 2u_m & \frac{2}{3}\alpha_1 \\
-2 u_m \alpha_1 & 2\alpha_1 & u_m
\end{pmatrix},
\end{align*}
and
\begin{align*}
\mathbf{S}(\mathbf{U}) = \cos(\theta)\Big(0,\,
h\tan(\theta) - \widetilde{\tau}_{xz}|_{\rm b} + \epsilon h \partial_x b,\,
-3\big( \widetilde{\tau}_{xz}|_{\rm b} + \mathcal{T}_1 \big)\Big)^T.
\end{align*}
The polynomial expansion of the flow velocity in the $\zeta$-coordinate is given by the linear function
\begin{equation}\label{eq:N=1u}
    \widetilde u(\zeta) = u_m + \alpha_1(1-2\zeta),
\end{equation}
where the only basis function is $\phi_1(\zeta) = 1-2\zeta$. From~\eqref{eq:N=1u}, we observe that $\partial_\zeta \widetilde{u}(\zeta) = -2\alpha_1$ for all $\zeta\in[0,1]$. In the case that $\alpha_1 < 0$, we have that $\sgn( \partial_\zeta \widetilde{u}(\zeta)) = 1$ and $\mathcal{J}^+ = (0,1)$ while $\mathcal{J}^-=\emptyset$. Note that condition $\widetilde{u}(0)=u_m + \alpha_1\geq 0$ will imply that $u_m \geq -\alpha_1$. Then, we can directly compute the bulk friction term \eqref{eq:bulkmuI} as
\begin{equation} \label{eq:integral:Tau1:N1}
\begin{split}
 \mathcal{T}_1 &= -4h \int_0^1 \xi^3 \left( \mu_{\rm s} + \frac{\mu_2 - \mu_{\rm s}}{1 + C_1 \xi} \right)d\xi \\
    &= -h\mu_{\rm s} - 4h(\mu_2-\mu_{\rm s}) \left( \frac{1}{3 C_1} -\frac{1}{2 C_1^2} + \frac{1}{C_1^3} -\frac{\log(1+C_1)}{C_1^4} \right),
\end{split}
\end{equation}
where
\begin{equation} \label{eq:tau_bot_allgr_N1}
C_1 = C_1(h,\alpha_1)
= \dfrac{\tilde{c}_I h^{3/2} }{2|\alpha_1|}.
\end{equation}
The bottom slip friction in this case is given by $\widetilde{\tau}_{xz}|_{\rm b} = (\nu_0/\lambda)(u_m + \alpha_1)$.
\begin{remark}
Similarly to this linear velocity profile case, we can compute algebraic expressions for $\mathcal{T}_i$ in the case $N=2$, for second order polynomial expansions of the velocity $\widetilde{u}$ in the transformed vertical coordinate $\zeta$. However, the computation of the integral terms becomes more involved due to more potential changes of the velocity sign and more cases need to be taken into account. In Appendix~\ref{sec:appendix:muI:N=2}, we have included some of the calculations related to the case $N=2$. The generalization of the computation of $\mathcal{T}_i$ for $i=1,\ldots,N$ for the case of $N \geq 3$ can be done in one of two ways: (1) at runtime, finding the roots of the velocity profile $\widetilde{u}$ and then splitting the integration into parts that can be integrated analytically or (2) using numerical integration such as Gaussian quadrature.
\end{remark}

\subsection{Equilibrium states analysis}\label{sec:equilibrium:stability}

One of the benefits of the modelling framework presented in this paper is the possibility to derive analytical equations even in the presence of complex rheologies. The derivation of the equations for the $\mu(I)$-rheology was performed explicitly in the previous section for $N=1$. Based on that, it is now possible to conduct further studies of the model structure. As one example, we consider the analysis of equilibrium states.

We thus illustrate the computation of equilibrium states of system (\ref{e:HSWME}) with the source term $\mathbf{S}$ given in (\ref{e:friction_terms}).
According to \cite{Yong1999}, equilibrium states are defined as states $\mathbf{U}$ at which the right-hand side source term $\mathbf{S}(\mathbf{U})$ vanishes. Letting $\mathbf{S}(\mathbf{U})=\boldsymbol{0}$, we have
\begin{align*}
    h\tan (\theta) = \widetilde{\tau}_{xz}|_{\rm b} - \widetilde{\tau}_{xz}|_{\rm s} + \epsilon h \partial_x b\qquad \text{and}\qquad \mathcal{T}_i =  (-1)^i \widetilde{\tau}_{xz}|_{\rm s} - \widetilde{\tau}_{xz}|_{\rm b},
\end{align*}
for $i=1,\dots,N$, which correspond to $N+1$ constraints for the vector of unknowns $\mathbf{U}=(h,hu_m,h\alpha_1,\dots,h\alpha_N)^{\rm T}$. Assuming flat topography $b=0$ and zero surface stress $\widetilde{\tau}_{xz}|_{\rm s}=0$ (as all the examples studied in this paper), the above equations reduce to
\begin{align} \label{eq:cond:equilibrium}
    h\tan (\theta) = \widetilde{\tau}_{xz}|_{\rm b}\qquad \text{and}\qquad \mathcal{T}_1= \mathcal{T}_2=\dots=\mathcal{T}_N= - \widetilde{\tau}_{xz}|_{\rm b}.
\end{align}

We now study the equilibrium states of the granular flow models with $\mu(I)$-rheology developed in Section \ref{sec:granular} for $N=1$. For the case of $N=1$ with $\alpha_1<0$, the equations in \eqref{eq:cond:equilibrium} combined with \eqref{eq:integral:Tau1:N1} give
\begin{align} \label{eq:aux:equilibrium}
\dfrac{\tan(\theta) - \mu_{\rm s}}{\mu_2-\mu_{\rm s}} =   4 \left( \frac{1}{3 C_1} -\frac{1}{2 C_1^2} + \frac{1}{C_1^3} -\frac{\log(1+C_1)}{C_1^4} \right) = 4\int_0^1 \dfrac{\varsigma^3}{1+C_1\varsigma} {\rm d}\varsigma,
\end{align}
where $C_1$ is defined in \eqref{eq:tau_bot_allgr_N1}. The integral in \eqref{eq:aux:equilibrium} can be bounded by
\begin{align*}
0<\dfrac{1}{1+C_1} < 4\int_0^1 \dfrac{\varsigma^3}{1+C_1\varsigma} {\rm d}\varsigma < 1,
\end{align*}
which allow us to conclude that \eqref{eq:aux:equilibrium} is valid if and only if
\begin{align*}
\mu_s < \tan(\theta) < \mu_2.
\end{align*}
Now, in the case of $\mu(I)$-rheology bottom friction \eqref{eq:muI:bottom}, we can determine all the equilibrium states. Due to the linear velocity profile with $N=1$, we evaluate $|\partial_\zeta\widetilde{u}(0)| = 2|\alpha_1|$ and ${\rm sgn}(\partial_\zeta\widetilde{u}(0)) = -1$ so that the friction term \eqref{eq:muI:bottom} reads
\begin{align*}
\widetilde{\tau}_{xz}|_{\rm b} & 
= -\bigg(\mu_{\rm s} + \frac{(\mu_2 - \mu_{\rm s})(2|\alpha_1|)}{ \tilde{c}_I h^{3/2} + 2|\alpha_1|}\bigg)h 
= -\bigg(\mu_{\rm s} + \frac{(\mu_2 - \mu_{\rm s})}{ 1+ C_1}\bigg)h.
\end{align*}
Then, from the second equation in~\eqref{eq:cond:equilibrium}, which reduces to $\mathcal{T}_1 = - \widetilde{\tau}_{xz}|_{\rm b}$, and taking into account the integral expression in \eqref{eq:aux:equilibrium}, we have 
\begin{align}\label{eq:aux:cond:Tau1}
  \dfrac{1}{1 + C_1} = \frac{4}{C_1} \left( \frac{1}{3} - \frac{1}{2C_1} + \frac{1}{C_1^2} - \frac{\log (1+C_1)}{C_1^3} \right) = 4\int_0^1 \dfrac{\varsigma^3}{1+C_1\varsigma}{\rm d}\varsigma > \dfrac{1}{1+C_1}.
\end{align}
As a consequence, the first equality in Equation \eqref{eq:aux:cond:Tau1} is only valid if $C_1\to \infty$, or equivalently if $\alpha_1 = 0$, see Equation \eqref{eq:tau_bot_allgr_N1}.
Therefore, the system has a non-trivial equilibrium if $h>0$, $u_m\ge 0$ and $\alpha_1=0$ if $\tan(\theta)=\mu_{\rm s}$, and does not have other equilibrium states with $h>0$.

The derived equilibrium state is thus characterized by a constant velocity profile $\widetilde u(\zeta) = u_m$ and only occurs under the condition that $\tan(\theta) = \mu_{\rm s} < \mu_2$. meaning that a flow with constant velocity is sliding down an inclined plane without any acceleration or change in velocity profile if the tangens of the inclination angle matches the coefficient $\mu_{\rm s}$. In that situation the bottom friction is precisely in equilibrium with the gravitational acceleration and the interior friction vanishes (due to zero velocity profile gradients). This demonstrates the consistency of the model with physical principles, enabled by the analytical derivation and the clear derivation framework.

\section{Numerical scheme} \label{sec:num:scheme}

The aim of this section is to develop a numerical method for the HSWME model that is suitable for simulations including wet-dry fronts. The main elements in this work are the variety of friction terms, which can become stiff in some cases and therefore need to be properly approximated.
For the spatial discretization over the domain $\Omega = (x_{\rm a},x_{\rm b})$, we consider here uniform meshes composed by $J$ cells $I_j = [x_{j-1/2}, x_{j+1/2}]$ of length $\Delta x=\tfrac{1}{J}( x_{\rm b} -x_{\rm a})$ with $x_{j\pm 1/2} = x_{\rm a} + (j\pm 1)\Delta x$ for all $j= 1,\dots,J$. On each cell and for all time $t>0$, we approximate the vector of unknowns by its cell average
\begin{align*}
 \mathbf{U}_j(t) = \dfrac{1}{\Delta x}\int_{x_{j-1/2}}^{x_{j+1/2}}\mathbf{U}(s,t)\,{\rm d}s,\qquad \text{for }j=1,\dots,J.
\end{align*}
Then, for the non-conservative product, we follow the approach in \cite{pares2006numerical}. At each cell interface $x_{j+1/2}$ with left $\mathbf{U}_j$ and right $\mathbf{U}_{j+1}$ cell values, we define the linearised Roe matrix $\mathbb{A}$ such that the following equality holds
\begin{align}\label{eq:linearized:Roe}
 \mathbb{A}\big(\mathbf{U}_j,\mathbf{U}_{j+1}\big) (\mathbf{U}_{j+1}-\mathbf{U}_{j}) = \int_{0}^{1} \mathbf{A}_{\rm H} \big(\mathbf{\Psi}(s; \mathbf{U}_j,\mathbf{U}_{j+1})\big) \dfrac{\partial \mathbf{\Psi}}{\partial s}\big(s; \mathbf{U}_j,\mathbf{U}_{j+1}\big)\,{\rm d}s,
\end{align}
where $\mathbf{\Psi}=\mathbf{\Psi}(s; \mathbf{U}_j,\mathbf{U}_{j+1})$ is a path connecting the two values $\mathbf{U}_j$ and $\mathbf{U}_{j+1}$. For the particular case of the linear path
\begin{align*}
 \mathbf{\Psi}(s;\mathbf{U}_j,\mathbf{U}_{j+1}) = \mathbf{U}_j + s(\mathbf{U}_{j+1} - \mathbf{U}_j),\qquad 0\leq s\leq 1,
\end{align*}
where $\partial_s \mathbf{\Psi} = \mathbf{U}_{j+1} - \mathbf{U}_j$, the linearised Roe matrix can be directly defined as
\begin{align}\label{eq:specific:Roe:matrixInt}
 \mathbb{A}(\mathbf{U}_j,\mathbf{U}_{j+1/2}) = \int_{0}^{1} \mathbf{A}_{\rm H} \big(\mathbf{U}_j + s(\mathbf{U}_{j+1} - \mathbf{U}_j)\big) \,{\rm d}s.
\end{align}
Taking into account the definition of $\mathbf{A}_{\rm H}$, where the nonlinear terms are products of the components of the vector $\mathbf{U}$, the composition $\mathbf{A}_{\rm H}(\mathbf{\Psi}(\,\cdot\,; \mathbf{U}_j,\mathbf{U}_{j+1}))$ becomes a second order polynomial in $s$. Therefore, the integral above can be either computed algebraically or using a quadrature rule exact for second order polynomials, such that
\begin{align}\label{eq:specific:Roe:matrix}
 \mathbb{A}_{j+1/2} = \sum_{l=1}^{n_{\rm G}} w_{l} \mathbf{A}_{\rm H} \big(\mathbf{U}_j + s_l(\mathbf{U}_{j+1} - \mathbf{U}_{j})\big),
\end{align}
where $n_{\rm G}$ is the number of quadrature points, $s_l$ and $w_l$ are the quadrature points and weights, respectively, for $l=1,\dots,n_{\rm G}$. Then, using a polynomial viscosity matrix (PVM) method \cite{castro2020well} for the non-conservative products, the semi-discrete approximation of the homogeneous PDE \eqref{e:HSWME}, is given by
\begin{equation}
\label{eq:semi_discrete}
\frac{{\rm d}\mathbf{U}_{j}}{{\rm d}t}= -\frac{1}{\Delta x}
\Big(
\mathbb{D}^{+}_{j-1/2} + \mathbb{D}^{-}_{j+1/2}\Big) + \mathbf{S}(\mathbf{U}_j),
\qquad \text{for }j=1,\dots,J,
\end{equation}
where $\mathbb{D}^{+}_{j-1/2}$ and $\mathbb{D}^{-}_{j+1/2}$ are the so-called fluctuations, which for the PVM method are defined as
\begin{align*}
 \mathbb{D}^{+}_{j-1/2} & = \mathbb{D}^{+}(\mathbf{U}_{j-1},\mathbf{U}_{j}):=\dfrac{1}{2}\Big\{\mathbb{A}_{j-1/2} + \mathbf{Q}_{j-1/2}\Big\}\big(\mathbf{U}_{j}-\mathbf{U}_{j-1}\big),\\[1ex]
 \mathbb{D}^{-}_{j+1/2} & = \mathbb{D}^{-}(\mathbf{U}_{j},\mathbf{U}_{j+1}):= \dfrac{1}{2}\Big\{\mathbb{A}_{j+1/2} - \mathbf{Q}_{j+1/2}\Big\}\big(\mathbf{U}_{j+1}-\mathbf{U}_{j}\big),
\end{align*}
with $\mathbf{Q}_{j+1/2}$ being the viscosity matrix, which can be defined as a polynomial function in terms of the matrix $\mathbb{A}_{j+1/2}$, that for the Price-C scheme \cite{Canestrelli2009} is defined as
\begin{align*}
 \mathbf{Q}_{j+1/2} = \tfrac{\Delta x}{2 \Delta t} \bB{I} + \tfrac{\Delta t}{2\Delta x}\mathbb{A}_{j+1/2}^2,
\end{align*}
for all $j = 0,1,\dots, J$, with $\bB{I}$ the identity matrix of size $N+2$. Ghost cells are supplemented to include boundary conditions, such that the numerical solution is extended to incorporate $\mathbf{U}_0$ and $\mathbf{U}_{J+1}$.

For the time approximation, we consider two discretizations: a fully explicit scheme given by forward Euler and a semi-implicit approach. For the explicit approximation, we have the following marching formula
\begin{align} \label{eq:explicit:in:time}
\mathbf{U}_j^{n+1} = \mathbf{U}_j^n - \dfrac{\Delta t}{\Delta x}\Big(\mathbb{D}^{+}(\mathbf{U}_{j-1}^n,\mathbf{U}_{j}^n) +
\mathbb{D}^{-}(\mathbf{U}_{j}^n,\mathbf{U}_{j+1}^n) \Big) +  \Delta t\, \mathbf{S}(\mathbf{U}_j^{n}).
\end{align}
For the semi-implicit variant, we keep the explicit approximation of the spatial derivatives and approximate the source term at $t^{n+1}$. Then, the marching formula for the semi-implicit scheme is given by
\begin{align*}
\mathbf{U}_j^{n+1} = \mathbf{U}_j^n - \dfrac{\Delta t}{\Delta x}\Big(\mathbb{D}^{+}(\mathbf{U}_{j-1}^n,\mathbf{U}_{j}^n) +
\mathbb{D}^{-}(\mathbf{U}_{j}^n,\mathbf{U}_{j+1}^n) \Big) +  \Delta t\, \mathbf{S}(\mathbf{U}_j^{n+1}),
\end{align*}
which is equivalent to solve Equation \eqref{e:HSWME} in two steps through the following splitting procedure:
\begin{align} \label{eq:semi-implicit:in:time}
\begin{aligned}
 \widecheck{\mathbf{U}}_j^{n+1} & = \mathbf{U}_j^n - \dfrac{\Delta t}{\Delta x}
 \Big(\mathbb{D}^{+}(\mathbf{U}_{j-1}^n,\mathbf{U}_{j}^n) +
\mathbb{D}^{-}(\mathbf{U}_{j}^n,\mathbf{U}_{j+1}^n)\Big),\\[1ex]
 \mathbf{U}_j^{n+1} & = \widecheck{\mathbf{U}}_j^{n+1} +  \Delta t\, \mathbf{S}(\mathbf{U}_j^{n+1}),
\end{aligned}
\end{align}
where $ \smash{\widecheck{\mathbf{U}}_j^{n+1}\in\mathbb{R}^{N+2}}$ corresponds to an intermediate variable used to compute the discrete unknown $\smash{\mathbf{U}_j^{n+1}}$. The solution of the second equation in \eqref{eq:semi-implicit:in:time} is computed using a Newton-Raphson solver. Although the explicit scheme does not require solving a non-linear system, the semi-implicit approach is more suitable for stiff source terms.
\smallskip

\noindent {\bf Wet-try treatment.} The so-called wet-dry fronts happen at interfaces between a wet cell ($h$ larger than $0$) and a dry cell ($h=0$). In practice, the case when $h$ tends to zero is problematic and needs a special treatment. This is mainly due to possible divisions by zero in the transformation between primitive and conservative variables and in the source term. From a physics perspective, when $h_j^n=0$ in a cell $I_j$, the only feasible solution is $\mathbf{U}_j^{n}=\mathbf{0}$. Then, we consider as a dry cell those in which $h_j^n\leq h_{\min}$ for a given minimum height $h_{\min}>0$, and to go from conservative to primitive variables, we use the following transformation \cite{Ciallella2025, Meister2016}:
\begin{align*}
(u_m)_j^n &= \dfrac{2 h_j^n}{(h_j^n)^2 + \max\{ (h_j^n)^2,h_{\rm min}\} }({\rm U}_2)_j^n,\\
(\alpha_i)_j^n &= \dfrac{2 h_j^n}{(h_j^n)^2 + \max\{ (h_j^n)^2,h_{\rm min}\}} ({\rm U}_{2+i})_j^n,
\end{align*}
for $i=1,\dots,N$, and set $(u_m)_j^n = (\alpha_i)_j^n = 0$ when $h_j^n<h_{\min}$. To deal with the source term, we use that the first component of $\mathbf{S}$ is in fact zero. This allow us to use the intermediate variable $\smash{\widecheck{h}_{j}^{n+1}}$ to determine the dry cells in which $\smash{\widecheck{h}_{j}^{n+1}}\leq h_{\min}$ such that we avoid computing the source term and set the mean velocity and moments to zero. This is particularly important when using the semi-implicit scheme \eqref{eq:semi-implicit:in:time}, in which we skip solving the Newton-Raphson solver in the dry cells.

For the time step size $\Delta t$, we compute the eigenvalues of $\mathbf{A}_{\rm H}(\mathbf{U}_j^n)$ using the formulas from \cite{Koellermeier2020} and use the formula $\Delta t\leq {\rm CFL}(\Delta x/\lambda_{\max})$ with $\lambda_{\max}$ the largest eigenvalue (in absolute value) of $\mathbf{A}_{\rm H}(\mathbf{U}_j^n)$, and CFL is a positive number.

\smallskip

\noindent {\bf Bulk friction approximation.}
The bulk friction $\mathcal{T}_i$ is defined in \eqref{eq:bulkmuI} through an integral over $0\leq \zeta\leq 1$, which for general expressions of $\widetilde{\tau}_{xz}$ may not be written in a closed algebraic form. This is the case in, for instance, the $\mu(I)$-rheology friction term in Section~\ref{sec:granular}. We approximate this integral via Gauss-Legendre quadrature as follows:
\begin{align} \label{eq:quadrature:TauI:muI}
 (\mathcal{T}_i)_j^n = h_j^n \sum_{l = 1}^{\tilde{n}_{\rm G}} \tilde{w}_l  (1-\tilde{s}_l)\widetilde{\tau}_{xz}(\tilde{s}_l; \mathbf{U}_j^n) \phi_i'(\tilde{s}_l), \qquad \text{for } i = 1,\dots,N,
\end{align}
where  $\tilde{w}_1,\dots,\tilde{w}_{\tilde{n}_G}$ and $\tilde{s}_1,\dots,\tilde{s}_{\tilde{n}_G}$ are the weights and nodes of a quadrature rule with $\tilde{n}_G$ points, respectively. The notation $\widetilde{\tau}_{xz}(\,\cdot\,; \mathbf{U}_j^n)$ is to emphasize that the friction term may depend on the unknown vector $\mathbf{U}_j^n$ which does not depend on $\zeta$. In the case of the $\mu(I)$-rheology, the bulk stress term $\widetilde{\tau}_{xz}|_{\rm bulk}$ can be written as
\begin{align*}
\widetilde{\tau}_{xz}|_{\rm bulk}(\tilde{s}_l;\mathbf{U}_j^n)  :=  \left( \mu_s + \frac{(\mu_2 - \mu_s) |(\bB{\alpha}_j^n)^{\rm T}\cdot \bB{\phi}'(\tilde{s}_l)|}{|(\bB{\alpha}_j^n)^{\rm T}\cdot \bB{\phi}'(\tilde{s}_l)| + C \sqrt{1-\tilde{s}_l}} \right) \sgn\big((\bB{\alpha}_j^n)^{\rm T}\cdot \bB{\phi}'(\tilde{s}_l)\big),
\end{align*}
where $\bB{\phi}=(\phi_1,\phi_2,\dots,\phi_N)^{\rm T}\in \mathbb{R}^N$. The shifted and scaled Legendre polynomials at the quadrature points $\bB{\phi}(\tilde{s}_l)$ for all $l=1,\dots,N$ can be precomputed before the time approximation is performed.

\section{Numerical simulations}\label{sec:simulations}

In this section, we present numerical simulations considering different scenarios for the friction models presented in the previous sections. The numerical examples are computed following the scheme described in Section~\ref{sec:num:scheme}, where all routines have been implemented in Fortran90, which can be found in \url{https://github.com/SWMEsolver/mainsolver}.

Through all the examples shown in this section, we use the following parameters
\begin{align*}
 &\epsilon = 0.01,\quad H = 0.1\,{\rm m},\quad L = 10\,{\rm m}, \quad U = \sqrt{gL}\approx 9.9045\,{\rm m/s},\quad \rho = 1200\,{\rm kg/m^3},\\
 &h_{\rm min} = 10^{-6},\quad J = 1000,
\end{align*}
with the exemption of Example 4, in which we use a different value for $\rho$.
In all examples, we use a rectangular block of material initially at rest as initial condition:
\begin{align} \label{eq:piecewise:h}
\begin{aligned}
&u_m(x,0) = \alpha_i(x,0) = 0,\qquad \text{ for } i = 1,\dots,N,\\
&h(x,0) = \begin{cases}
           0.08 & \text{ if}\quad 0.3 \leq x \leq 0.5,\\
           h_{\min} & \text{otherwise}.
          \end{cases}
\end{aligned}
\end{align}
For the boundary conditions, we use transmissive boundary conditions at the left and right boundary of the domain, meaning that the value at each ghost cell is equal to that at the respective limit, that is
\begin{align*}
 \mathbf{U}_0^n = \mathbf{U}_1^{n},\qquad \mathbf{U}_{J+1}^n = \mathbf{U}_{J}^n, \qquad \text{for }n\geq 0.
\end{align*}
For the path integral in \eqref{eq:specific:Roe:matrix}, we consider a linear path in all examples, and use the following third order Gauss-Legendre quadrature rule for its approximation:
\begin{align*}
w_1 &= \tfrac{5}{18},\quad  w_2 = \tfrac{8}{18},\quad  w_3 = \tfrac{5}{18},\\
s_1 &= \tfrac{1}{2}\Big( 1 -\tfrac{\sqrt{15}}{5}\Big),\quad s_2 = \tfrac{1}{2},\quad s_3 = \tfrac{1}{2}\Big( 1 + \tfrac{\sqrt{15}}{5}\Big),
\end{align*}
where $n_{\rm G} = 3$. In the first two examples, we use the semi-implicit scheme, whereas in the remaining examples we employ the fully explicit one. The CFL constant used for the first three examples is ${\rm CFL} = 0.05$, and ${\rm CFL} = 0.01$ for the fourth example. In the case of the semi-implicit scheme, we implemented the Newton-Raphson solver making use of central differences to approximate the respective Jacobian matrix. The tolerance for the nonlinear solver is set to $\texttt{tol} = 10^{-6}$.

\begin{figure}[t]
 \includegraphics[scale=1]{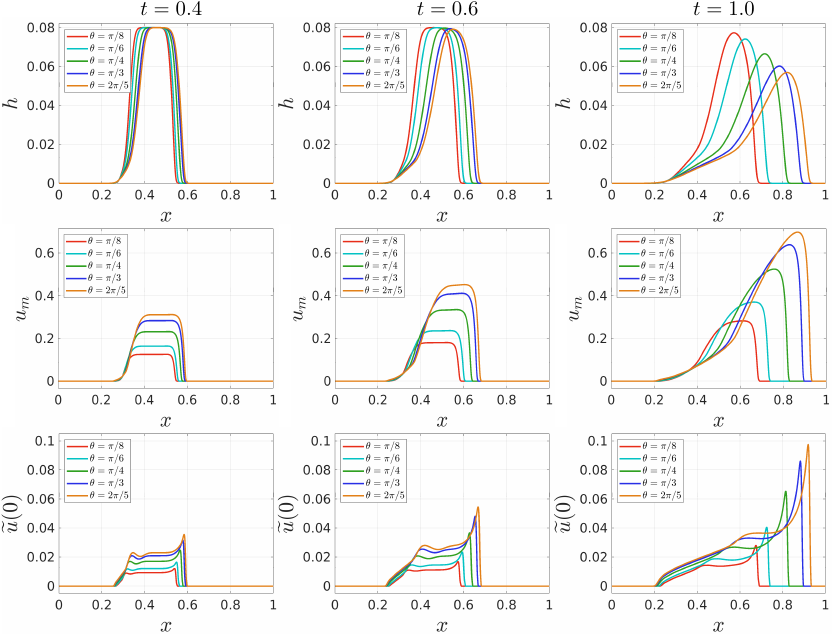}
  \caption{Example 1: Inclination effect on the height $h$ (top row), average velocity $u_m$ (middle row) and bottom velocity $\widetilde{u}(0)$ (bottom row) for the Newtonian-slip friction at $t=0.4$ (left column), $t=0.6$ (middle column) and $t=1$ (right column). In all the examples $\Lambda = 10^{-4}$ and $\eta=0.01$.}\label{fig:Figure1}
\end{figure}

\subsection{Example 1: Inclination angle effects of Newtonian slip model}\label{sec:example1}

To show the effect that the inclination angle $\theta$ has on the flow, we simulate a sliding avalanche starting with the piecewise constant height $h(x,0)$ given by the initial condition \eqref{eq:piecewise:h} initiating from rest, i.e. zero initial velocity. We consider the case of a Newtonian-slip friction described in Section~\ref{sec:Newtonian_Slip_Flow} with fixed slip length $\Lambda = 10^{-4}$ and viscosity constant $\eta = 0.01$ for the case $N=2$.

In Figure~\ref{fig:Figure1}, we show snapshots of the numerically computed height $h$, average velocity $u_m$, and bottom velocity $\widetilde{u}(0)$ (see \eqref{eq:bottom:velocity}) at three different times, for different angles from $\theta=\pi/8$ to $\theta = 2\pi/5$. As expected, the inclination angle increases the horizontal velocity and therefore $h$ evolves faster. We also observe that there is a clear front in the bottom velocity $\widetilde{u}(0)$, which jumps from zero to a peak value that decreases from right to left in the $x$-axis. This peak in the bottom velocity originates from the initial collapse of the material under its own weight. The average velocity $u_m$, on the other hand, exhibit a smoother behaviour for the different angles simulated, and no sharp peaks are observed. In all examples, the velocity remains positive and the small non-zero values to the left of $x=0.3$ in the height can be explained due to numerical errors.

\begin{figure}[t]
 \includegraphics[scale=1.07]{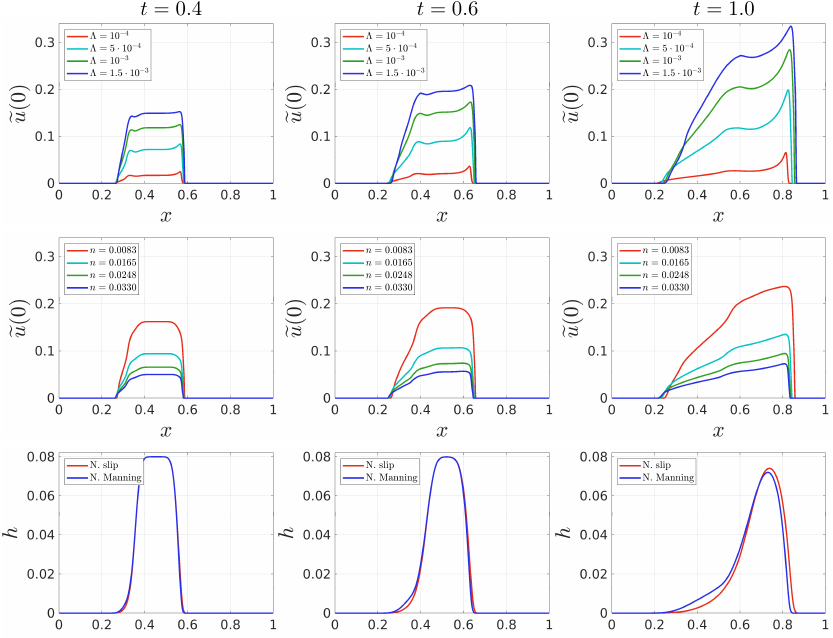}
 \caption{Example 2: Bottom friction effect on the bottom velocity for the Newtonian-slip (top row) and Newtonian-Manning (middle row) for $\theta = \pi/4$ and different values of $\Lambda$ and $n$ at $t=0.4$ (left column), $t=0.6$ (middle column) and $t=1$ (right column). Bottom row: comparison of the height $h$ between Newtonian-slip with $\Lambda=0.0015$ and Newtonian-Manning with $n = 0.0165$.}\label{fig:Figure2}
\end{figure}

\subsection{Example 2: Bottom friction effects of Newtonian Manning model}\label{sec:example2}

To compare the effect of different types of bottom friction $\widetilde{\tau}_{xz}|_{\rm b}$ for a common bulk friction, we consider the case of a Newtonian bulk friction for slip (Section~\ref{sec:Newtonian_Slip_Flow}) and Manning (Section~\ref{sec:Newtonian_Manning}) bottom friction. To do so, we employ different values of the slip length $\Lambda$ and the Manning coefficient $n$ and determine their effect on the numerically computed bottom velocity $\widetilde{u}(0)$ for the case $N=2$.

In Figure~\ref{fig:Figure2}, we compare the bottom velocity $\widetilde{u}(0)$ for the Newtonian-slip and Newtonian-Manning models with the respective parameters:
\begin{align*}
\Lambda & = 0.0001, 0.0005, 0.0010, 0.0015,\\
n & = 0.0083, 0.0165, 0.0248, 0.0330,
\end{align*}
at three time points $t=0.4$, $t=0.6$ and $t=1$. Both types of friction models exhibit qualitative differences in terms of the variation of the bottom velocity $\widetilde{u}(0)$. While the Newtonian-slip shows distinctive peaks at the leading front, the Manning law results in smoother solutions and the magnitude of the velocity is clearly smaller. In the third row of Figure~\ref{fig:Figure2}, we show that choosing $\Lambda=1.5\cdot 10^{-3}$ for the Newtonian-slip and $n=0.0165$ for the Newtonian-Manning, the height evolves similarly in both cases up to a certain time, in which the difference gets more pronounced due to the faster movement in the Newtonian-slip.

\begin{figure}
 \includegraphics[scale=1.07]{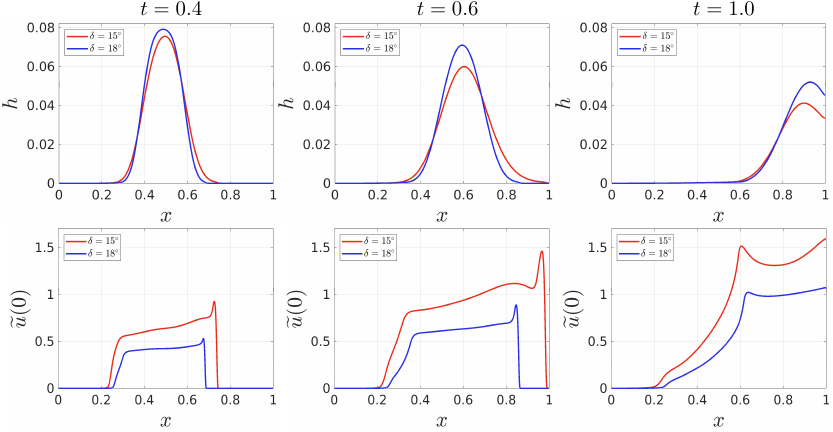}
\caption{Example 3: Comparison of the height $h$ (top row) and bottom velocity $\widetilde{u}(0)$ (bottom row) computed with the Savage-Hutter friction model for two different bed friction angles $\delta=15^\circ$ (red) and $\delta = 18^\circ$ (blue) at $t=0.4$ (left column), $t=0.6$ (middle column) and $t=1$ (right column). The inclination angle used is $\theta = \pi/4$ and inner friction angle $\varphi = 20^\circ$.}\label{fig:Figure3}
\end{figure}

\subsection{Example 3: Bed friction effects of Savage-Hutter model}\label{sec:example3}

Before we discuss the simulation results for this example, a brief discussion of the assumption of positive bottom velocity $\widetilde{u}(0)$ is in order: The Savage-Hutter model in Section~\ref{sec:Savage-Hutter_model} differs from Newtonian based models mainly because the bottom friction depends on the bed friction angle $\delta$, while the bulk friction is affected by the inner angle $\varphi$. Both angles are such that $\delta < \theta$ and $\delta \leq \varphi$. In the event of sliding avalanches, or the case of the initial condition given by \eqref{eq:piecewise:h}, the velocity is expected to be positive and the bulk friction is assumed to be given by the expression in \eqref{eq:SHbulk}. Conversely, one may expect that given a bulk friction term $\mathcal{T}_i$ in the form of \eqref{eq:SHbulk}, the velocity should be positive. However, the latter is not always true and we may ensure that the inner angle $\varphi$ is close enough to $\delta$ to verify this condition. This can be explained looking at the source term for the moments, which contains the difference $ \tan(\varphi)-\tan(\delta)$ that for larger values of $\varphi$ affects the moments such that the velocity could become negative.

In Figure~\ref{fig:Figure3}, we present the comparison of the height $h$ and bottom velocity $\widetilde{u}(0)$ for the Savage-Hutter friction model at two different bed friction angles $\delta = 15^\circ$ and $\delta = 20^\circ$ for a common inner angle $\varphi = 20^\circ$ and $\theta = \pi/4$ and $N=2$. We observe that the bed friction angle $\delta$ has a strong impact on the solution with smaller values of $\delta$ implying larger velocity profiles. We also observe that the deformation of the initial condition is more pronounced than for the Newtonian examples (Sections \ref{sec:example1} and \ref{sec:example2}), meaning that even at early times, the height exhibits a Gaussian shape. For the parameters chosen, the bottom velocity $\widetilde{u}(0)$ tends to be higher than in previous examples and already at $t=1$ material has reached the right boundary. A narrow peak is also observed in the bottom velocity near the leading front like those observed in the Newtonian-slip case. Again, this could be attributed to the quick initial collapse of the material under its own weight, which then propagates downhill.

\subsection{Example 4: Varying moments and runoff of $\mu(I)$-rheology model} \label{sec:example4}

In this example, we consider a bottom slip law and the definition of the $\mu(I)$-rheology based bulk friction term in \eqref{eq:bulkmuI} without making an assumption on the sign of the velocity profile. Therefore the bulk friction $\mathcal{T}_i$ needs to be computed using a quadrature rule as described in \eqref{eq:quadrature:TauI:muI}. To ensure that the approximation of this integral is sufficiently accurate and will not affect the numerical approximation, we use eight order Gauss-Legendre quadrature rule, that is $\tilde{n}_{\rm G} = 8$, with quadrature points
\begin{align*}
&
\tilde{s}_1 = 0.01985,\quad
\tilde{s}_2 = 0.10167,\quad
\tilde{s}_3 = 0.23723,\quad
\tilde{s}_4 = 0.40828,\\
&
\tilde{s}_5 = 0.59172,\quad
\tilde{s}_6 = 0.76277,\quad
\tilde{s}_7 = 0.89833,\quad
\tilde{s}_8 = 0.98015,
\end{align*}
and quadrature weights
\begin{align*}
&\tilde{w}_1 = 0.05061,\quad
\tilde{w}_2 = 0.11119,\quad
\tilde{w}_3 = 0.15685,\quad
\tilde{w}_4 = 0.18134,\\
&
\tilde{w}_5 = 0.18134,\quad
\tilde{w}_6 = 0.15685,\quad
\tilde{w}_7 = 0.11119,\quad
\tilde{w}_8 = 0.05061.
\end{align*}
Further parameters used in this example are $I_0 = 0.279$, $d_s = 0.7\,{\rm mm}$, $\mu_s  = 0.48$, $\mu_2  = 0.73$, $\rho = 1550\,{\rm kg/m^3}$ and $\rho_s = 2500\,{\rm kg/m^3}$ (see \cite{FernandezNieto2016}), and for the bottom friction we set $\Lambda = 0.001$ and $\eta_0=
0.001$ (see \cite{KoellermeierQian2022}).

\begin{figure}[t]
 \includegraphics[scale=1.07]{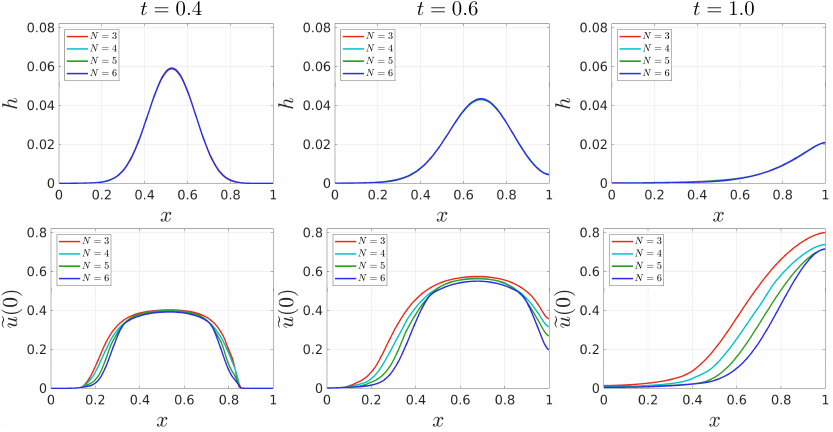}
\caption{Example 4: Comparison of the height $h$ (top row) and bottom velocity $\widetilde{u}(0)$ (bottom row) computed with the $\mu(I)$-rheology friction model with bottom slip law for different number of moments $N=3,4,5,6$ at $t=0.4$ (left column), $t=0.6$ (middle column) and $t=1$ (right column). The inclination angle used is $\theta = \pi/4$.}\label{fig:Figure4}
\end{figure}

We first investigate the variations of the numerical solutions with respect to the number of moments $N$. In Figure~\ref{fig:Figure4}, we show snapshots of the simulated granular flow using the $\mu(I)$-rheology based friction for $N=3,4,5,6$ at three time points $t=0.4, 0.6, 1.0$. As also observed in tests with different number of moments in \cite{kowalski2018moment}, the height $h$ (first row), is almost not changing, while the velocity or in this case the bottom velocity $\widetilde{u}(0)$ (bottom row) is still varying as the number of moments increases. It appears that the changes between the respective models decrease, which indicates convergence.

\begin{figure}[t]
\includegraphics[scale=0.9]{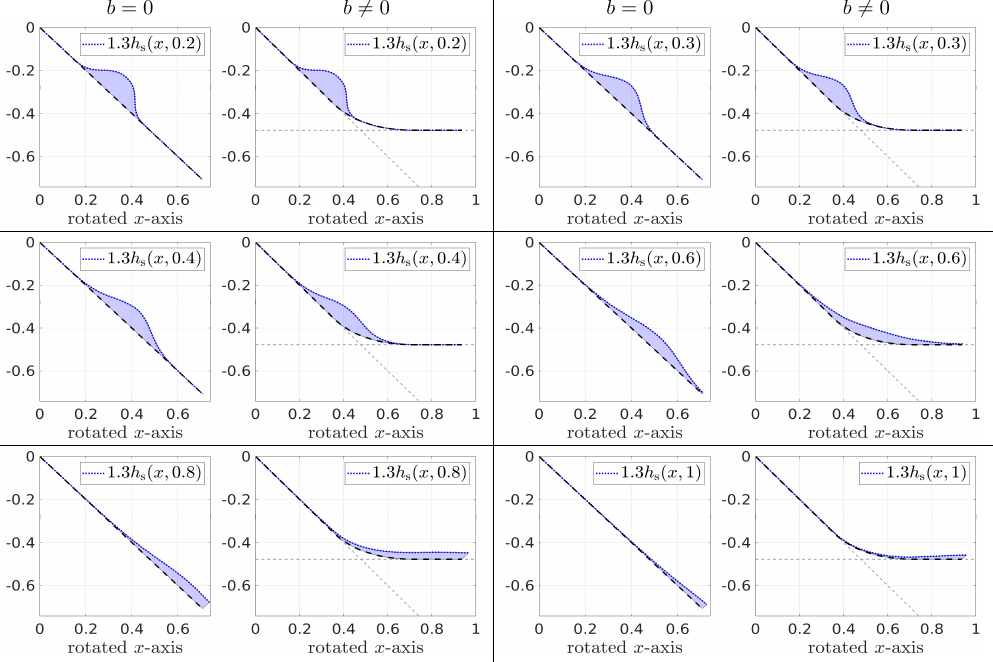}
\caption{Example 4: Free surface $h_{\rm s}$ computed with the $\mu(I)$-rheology friction model  with bottom slip law with $b = 0$ and $b\neq 0$ defined by \eqref{eq:bottom:curve} at different times from $t = 0.2$ to $t=1$. The inclination angle used is $\theta = \pi/4$. The plots are organised by pairs, and the time points increase from left to right and top to bottom.}\label{fig:Figure5}
\end{figure}

\begin{figure}[t]
\includegraphics[scale=1]{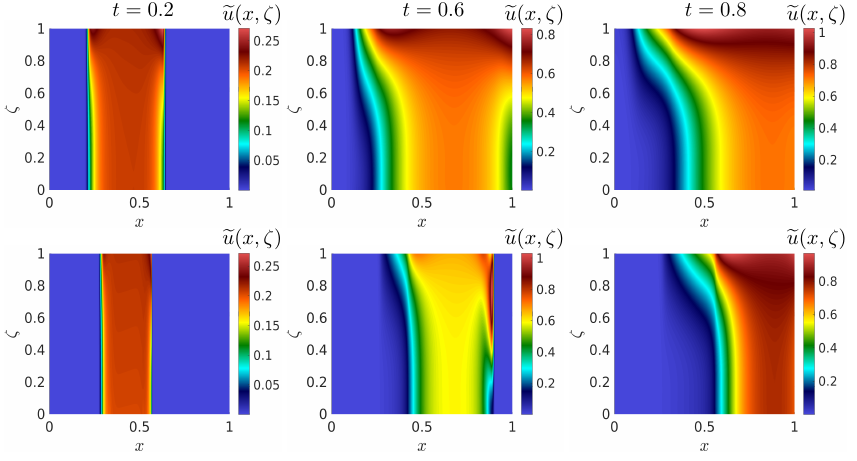}
\caption{Example 4: Velocity $\widetilde{u}$ computed with the $\mu(I)$-rheology friction model  with bottom slip law with $b = 0$ (top row) and $b\neq 0$ defined by \eqref{eq:bottom:curve} (bottom row) at $t = 0.2$ (left column), $t=0.6$ (middle column) and $t = 0.8$ (right column). The inclination angle used is $\theta = \pi/4$. The colormap is scaled with respect to the maximum value of $\widetilde{u}$ at each time.}\label{fig:Figure6}
\end{figure}

Furthermore, to test the influence of a variable topography in our model, we consider the following bottom curve which models the runoff from an inclined plane into a vertical plane:
\begin{align} \label{eq:bottom:curve}
 b(x) = \begin{cases}
            0 &\text{ if }  x < 0.5,\\
            \tfrac{10}{7}\tan(\theta) (x - 0.5)^2  &\text{ if } 0.5 \leq x \leq 0.85,\\
            \tan(\theta) (x - 0.675) &\text{ if } x > 0.85,
        \end{cases}
\end{align}
and compute $\partial_x b_{j}\approx \tfrac{1}{2}(\partial_x b(x_{j-1/2}) + \partial_x b(x_{j+1/2}))$ for the approximation of the source term $\mathbf{S}$. In Figure~\ref{fig:Figure5}, we present the simulation of the surface height $h_s = h+b$ for the case of flat topography $b=0$ and for variable bathymetry in a rotated $x$-axis for the case of $N=3$. The snapshots from $t=0.2$ to $t=1$ show the evolution of the flow surface amplified by a factor $1.3$ to improve their visualization. The velocity $\widetilde{u}$, which is given by the polynomial expansion \eqref{eq:mom_ansatz}, is plotted in Figure~\ref{fig:Figure6} as a function of $x$ and $\zeta$. We observe that the bottom topography directly affects the velocity of the granular flow, slowing it down as the flow approaches the plateau for $x>0.85$. Figure~\ref{fig:Figure6} also shows that in both cases the velocity varies with respect to $\zeta$ and that in general this function is increasing from bottom to top.

\section{Conclusions}\label{sec:conclusions}

In this paper we extend the shallow water moment model to granular flows on inclined planes with general friction laws via a general derivation framework.
A structured modelling procedure is developed to include arbitrary types of friction models. We study in particular the case of the Savage-Hutter friction and $\mu(I)$-rheology based granular flows. For the Savage-Hutter friction model, which is defined via the bed and inner friction angles, we observe that an additional condition on the inner friction needs to be imposed. For large inner friction, unphysical negative velocities may arise due to the influence of this angle in the moment equations. The general $\mu(I)$-rheology is included in the model by approximating the integral related to the bulk friction via a quadrature rule. Regarding the bottom friction, our flexible approach allows to make different choices. In one example, we combine the $\mu(I)$-rheology bulk friction with a slip condition for the bottom. However, one can also use a Coulomb friction law for the bottom, but the influence of the bed friction angle $\delta$ needs to be analysed.

This work opens many possible directions for future research. Further studies can be conducted to determine the limiting values of the parameters bottom friction angle $\delta$ and inner friction angle $\varphi$ for the Savage-Hutter model (or $\mu$ in the case of the Coulomb-type friction) such that the velocity in Section~\ref{sec:Savage-Hutter_model} and increasing behaviour in Section~\ref{sec:Coulomb-type_const} are fulfilled. In terms of numerical approximation, different PVM methods can be implemented and tested, and higher order schemes can also open a pathway for extension. In particular, a proper derivation of a well-balanced scheme to preserve steady state solutions in the case of non-zero bottom topography, as done in \cite{FernandezNieto2018}, is of research interest.

%\section*{Author contributions}
%\textbf{Julio Careaga:} Conceptualization, Methodology, Investigation, Formal Analysis, Software, Validation, Visualization, Writing---original draft, Writing---review and editing.
%\textbf{Qian Huang:} Conceptualization, Methodology, Investigation, Formal Analysis, Writing---original draft,  Writing---review and editing.
%\textbf{Julian Koellermeier:} Funding Acquisition, Conceptualization, Methodology, Investigation, Formal Analysis, Writing---original draft, Writing---review and editing.

\section*{Acknowledgments}
JK and JC are supported by the Dutch Research Council (NWO) through the ENW Vidi project HiWAVE with file number VI.Vidi.233.066. QH is supported by the Deutsche Forschungsgemeinschaft (DFG, German Research Foundation) - SPP 2410 \textit{Hyperbolic Balance Laws in Fluid Mechanics: Complexity, Scales, Randomness} (CoScaRa).

\appendix
\section{Appendix}
\subsection{Derivation of Equation \eqref{e:dyn_BCb_nd}} \label{sec:stress_tensor_bottom}

We first recall the definition of $\bB{\tau}$ in \eqref{eq:tau:components} and the normal and tangent vectors in \eqref{eq:normal:tangent:b} and \eqref{eq:normal:tangent:s}, respectively. Then, it is straightforward to observe that:
\begin{align} \label{eq:identity:pI:tan}
\big((\bB{\tau}|_{\rm b} - p|_{\rm b}\bB{I})\bB{n}_{\rm b}\big)_{\rm tan} &= \big(\bB{\tau}|_{\rm b}\bB{n}_{\rm b}\big)_{\rm tan}, \qquad
\big((\bB{\tau}|_{\rm s} - p|_{\rm s}\bB{I})\bB{n}_{\rm s}\big)_{\rm tan} = \big(\bB{\tau}|_{\rm s}\bB{n}_{\rm s}\big)_{\rm tan}.
\end{align}
On the other hand, multiplying \eqref{e:dyn_BC_bottom} by the tangential vector at the bottom bathymetry $\bB{t}_{\rm b}$, and analogously multiplying \eqref{e:dyn_BC_surface} by $\bB{t}_{\rm s}$, and using \eqref{eq:identity:pI:tan}, we have
\begin{align*}
 T_{\rm b} &= (\bB{\tau}|_{\rm b}\bB{n}_{\rm b})_{\rm tan}\cdot \bB{t}_{\rm b} = (\bB{\tau}|_{\rm b}\bB{n}_{\rm b})\cdot \bB{t}_{\rm b},\qquad
 T_{\rm s}  =  (\bB{\tau}|_{\rm s}\bB{n}_{\rm s})_{\rm tan}\cdot \bB{t}_{\rm s} = (\bB{\tau}|_{\rm s}\bB{n}_{\rm s})\cdot \bB{t}_{\rm s}.
\end{align*}
Then, from the definition of $\bB{n}_{\rm b}$ and $\bB{t}_{\rm b}$ in \eqref{eq:normal:tangent:b}, we obtain
\begin{align*}
T_{\rm b}
& = \dfrac{1}{(\partial_x b)^2 + 1} \begin{pmatrix}
\tau_{xx}|_{\rm b}\partial_x b - \tau_{xz}|_{\rm b} \\[1ex]
\tau_{xz}|_{\rm b}\partial_x b -  \tau_{zz}|_{\rm b}
\end{pmatrix}\cdot
\begin{pmatrix}
    1\\
     \partial_x b
\end{pmatrix}
\\
& = \dfrac{1}{(\partial_x b)^2 + 1} \Big(\tau_{xx}|_{\rm b}\partial_x b - \tau_{xz}|_{\rm b}  + \partial_x b(\tau_{xz}|_{\rm b}\partial_x b -  \tau_{zz}|_{\rm b})\Big)\\
& = \dfrac{1}{(\partial_x b)^2 + 1} \Big(\big(\tau_{xx}|_{\rm b}-\tau_{zz}|_{\rm b}\big)\partial_x b - \big((\partial_x b)^2-1\big) \tau_{xz}|_{\rm b} \Big).
\end{align*}
Similarly, we can find the following expression for the tangential component at the surface
\begin{align*}
T_{\rm s}
& = -\dfrac{1}{(\partial_x h_{\rm s})^2 + 1} \Big(\big(\tau_{xx}|_{\rm s}-\tau_{zz}|_{\rm s}\big)\partial_x h_{\rm s} - \big((\partial_x h_{\rm s})^2-1\big) \tau_{xz}|_{\rm s} \Big).
\end{align*}
Following the scaling \eqref{eq:nondim} in Section~\ref{sec:non-dim:model}, we can conclude that $\partial_{x} b$ and $ \partial_{x} h_{\rm s}$ are of order $\mathcal{O}(\epsilon)$ while the products $\big(\tau_{xx}|_{\rm b}-\tau_{zz}|_{\rm b}\big)\partial_x b$ and $\big(\tau_{xx}|_{\rm s}-\tau_{zz}|_{\rm s}\big)\partial_x h_{\rm s}$ are of order $\mathcal{O}(\epsilon^2)$, therefore after neglecting terms of order $\mathcal{O}(\epsilon^2)$, the dimensionless versions of $T_{\rm b}$ and $T_{\rm s}$ are
\begin{align*}
 T_{\rm b}^* = \big(\rho g \cos(\theta)H\big) \tau_{xz}^*|_{\rm b},\qquad  T_{\rm s}^* = -\big(\rho g \cos(\theta)H\big) \tau_{xz}^*|_{\rm s}.
\end{align*}

\subsection{Case $N=2$, $\mu(I)$-rheology friction} \label{sec:appendix:muI:N=2}

We study here the $\mu(I)$-rheology friction model in Section~\ref{sec:granular} for $N=2$, i.e. a quadratic expansions of $\widetilde{u}$ in the transformed vertical coordinate $\zeta$.
In this case, the moments are $\alpha_1$ and $\alpha_2$, and the vector of unknowns is defined as $\mathbf{U} = (h,hu_m,h\alpha_1,h\alpha_2)^{\rm T}$. Furthermore, the system matrix $\mathbf{A}_{\rm H}$ and source term $\mathbf{S}$ in \eqref{e:HSWME} are given by
\begin{equation*}
\mathbf{A}_{\rm H}(\mathbf{U}) =
        \left(
        \begin{array}{ccccccc}
         0 & 1 & 0 & 0  \\
        \epsilon \cos(\theta) h-u_m^2-\frac{1}{3} \alpha_1^2& 2 u_m & \frac{2}{3} \alpha_1 & 0  \\
         -2u_m  \alpha_1 & 2 \alpha_1 & u_m & \frac{3}{5} \alpha_1  \\
         -\frac{2}{3} \alpha_1^2 & 0 & \frac{1}{3} \alpha_1 & u_m
        \end{array}
        \right),
\end{equation*}
and
\begin{align*}
 \mathbf{S}(\mathbf{U}) =
 \cos(\theta)\Big(0,\,
  h\tan(\theta)- \widetilde{\tau}_{xz}|_{\rm b}  + \epsilon h \partial_x b,\,
-3\big( \widetilde{\tau}_{xz}|_{\rm b} + \mathcal{T}_1 \big),
-5\big( \widetilde{\tau}_{xz}|_{\rm b} + \mathcal{T}_2 \big)\Big)^{\rm T}.
\end{align*}
The two polynomials that form the basis in this case are $\phi_1(\zeta) = 1-2\zeta$ and $\phi_2(\zeta) = 1-6\zeta+6\zeta^2$, and the polynomial expansion \eqref{eq:mom_ansatz} reduces to:
\begin{align*}
    \widetilde u(\zeta) 
    %&= u_m + \alpha_1 \phi_1(\zeta) + \alpha_2 \phi_2(\zeta) 
    = 6\alpha_2\zeta^2 -2(3\alpha_2 + \alpha_1)\zeta +\alpha_2+\alpha_1+ u_m,
\end{align*}
Then, we can readily compute the derivative $\partial_\zeta \widetilde u(\zeta) = 2(6\alpha_2\zeta - \alpha_1 -3\alpha_2)$, and observe that $\widetilde u(\zeta)$ attains its extreme value at
\begin{equation}\label{eq:extreme}
    \zeta^*=\frac{1}{2}\Big(1+\frac{\alpha_1}{3\alpha_2}\Big),
\end{equation}
with which we can write $\partial_\zeta \widetilde u(\zeta) =  12\alpha_2(\zeta-\zeta^*)$. One may additionally require $\widetilde u(\zeta)>0$ on $[0,1]$, or equivalently $u(0) = u_m+\alpha_1+\alpha_2>0$ and $u(1) = u_m-\alpha_1+\alpha_2>0$, which can be reduced to $u_m > |\alpha_1| - \alpha_2$. 

Now, we can identify sub-cases depending on the sign of $\partial_{\zeta}\widetilde{u}$, and whether $\zeta^*$ lies in the interval $(0,1)$ or not. Condition $\partial_\zeta \widetilde{u}>0$ holds in two scenarios: $\alpha_2>0$ and $\zeta>\zeta^*$, or $\alpha_2<0$ and $\zeta<\zeta^*$. In what follows, we study the sub-cases depending on $\zeta^*$ for which $\partial_\zeta \widetilde{u}>0$.

\bigskip

\noindent\textit{Case I ($\zeta^* \leq 0$ or $1\leq \zeta^* $):} We observe that the velocity satisfies condition $\partial_\zeta \widetilde{u}>0$ in the interval $(0,1)$, in the only possible sub-cases:  
\begin{align*}
( \alpha_2>0 \text{ and }\zeta^*\leq 0)\quad \text{or}\quad (\alpha_2 < 0 \text{ and } 1\leq \zeta^*).
\end{align*}
Then, we can compute the two additional bulk friction terms \eqref{eq:bulkmuI} through the formula \eqref{eq:bulkmuI_N} to obtain
\begin{align}
    \mathcal{T}_1  &= -h\mu_s - 4h(\mu_2-\mu_s) \int_0^1 \frac{A\xi^5+B\xi^3}{A\xi^2-C\xi+B} {\rm d}\xi,\label{eq:Tau1}\\[1ex]
\mathcal{T}_2      &= -h\mu_s - 12h(\mu_2-\mu_s) \int_0^1 \frac{(2\xi^5-\xi^3)(A\xi^2+B)}{A\xi^2-C\xi+B} {\rm d}\xi, \label{eq:Tau2}
\end{align}
where $A=12\alpha_2$, $B=2\alpha_1-6\alpha_2$, and $C = \tilde{c}_I h^{3/2}$. We observe that the remaining integral terms in \eqref{eq:Tau1} and \eqref{eq:Tau2} are the integrals of rational functions in $\xi$, which can be computed making use of partial fraction decomposition. For instance, for the integral in \eqref{eq:Tau1}, we have
\begin{align*}
\begin{aligned}
    \int_0^1 \frac{A\xi^5+B\xi^3}{A\xi^2-C\xi+B} {\rm d}\xi
    &= \frac{1}{4} + \frac{C}{3A} + \frac{C^2}{2A^2} - \frac{C(AB-C^2)}{A^3} + \frac{C}{A^3} \int_0^1 \frac{\omega_1\xi+\omega_2}{A\xi^2-C\xi+B} {\rm d}\xi,
\end{aligned}
\end{align*}
where $\omega_1=C^3 - 2ABC$ and $\omega_2=AB^2-BC^2$. In turn, the last integral on the right-hand side of the above equality can be computed as follows
\begin{align*}
\begin{aligned}
    \int_0^1 &\frac{\omega_1\xi+\omega_2}{A\xi^2-C\xi+B} {\rm d}\xi
    = \frac{\omega_1}{2A} \ln \left| \frac{A-C+B}{B} \right| +
    \begin{cases}
     \chi_1 \text{artanh}(\chi_2) & \text{if } C^2 \neq 4AB, \\[1ex]
     \dfrac{4A\omega_2+2C\omega_1}{C(C-2A)} & \text{if } C^2=4AB,
    \end{cases}
\end{aligned}
\end{align*}
where the constants $\chi_1$ and $\chi_2$ are given by
\begin{align*}
\chi_1 = \frac{2A\omega_2 + C\omega_1}{A\sqrt{|C^2-4AB|}}, \qquad  \chi_2 = \frac{\sqrt{|C^2-4AB|}}{2B-C}.
\end{align*}
The integral term in \eqref{eq:Tau2} can be treated analogously.

\bigskip

\noindent\textit{Case II ($0 < \zeta^* <1$):} We first observe that from the definition of $\zeta^*$ in \eqref{eq:extreme}, we get the bound $|\alpha_1/\alpha_2|\leq 3$. Now, the sign of $\partial_\zeta \widetilde{u}(\zeta)$ will change within the interval $(0,1)$ depending on the sign of $\alpha_2$. When $\alpha_2<0$, we have that $|\alpha_1|\leq -3\alpha_2$, and  
\begin{align} \label{eq:dudzeta:caseII}
\partial_\zeta\widetilde{u}(\zeta)=12\alpha_2(\zeta-\zeta^*)
\begin{cases}
>0 & \text{if }\zeta< \zeta^*,\\
\leq 0 & \text{if }\zeta\geq \zeta^*,
 \end{cases}
\end{align}
from which we can notice that the intervals of integration in \eqref{eq:bulkmuI_N} are $\mathcal{J}^+ = (0,\zeta^*)$ and $\mathcal{J}^{-}=(\zeta^*,1)$. Then, the additional bulk friction term \eqref{eq:bulkmuI}, or equivalently \eqref{eq:bulkmuI_N}, can be written as
\begin{align}
&\begin{aligned}
    \mathcal{T}_1 = &-4h
\int_{\xi^*}^1 \xi^3 \left( \mu_s + \frac{(\mu_2-\mu_s)(A\xi^2+B)}{A\xi^2-C\xi+B} \right) {\rm d}\xi  \\
& +4h \int_0^{\xi^*} \xi^3 \left( \mu_s + \frac{(\mu_2-\mu_s)(A\xi^2+B)}{A\xi^2+C\xi+B} \right) {\rm d}\xi,
\end{aligned}\\[1ex]
&\begin{aligned}
\mathcal{T}_2 = &-12h
\int_{\xi^*}^1 \xi^3 (2\xi^2-1) \left( \mu_s + \frac{(\mu_2-\mu_s)(A\xi^2+B)}{A\xi^2-C\xi+B} \right) {\rm d}\xi \\
& +12h \int_0^{\xi^*} \xi^3 (2\xi^2-1) \left( \mu_s + \frac{(\mu_2-\mu_s)(A\xi^2+B)}{A\xi^2+C\xi+B} \right) {\rm d}\xi,
\end{aligned}
\end{align}
where $\xi^* = (1-\zeta^*)^{1/2}$. These integrals can be evaluated either analytically or numerically through Gauss-Legendre quadrature rule.
The case $\alpha_2> 0$ can be studied analogously exchanging the inequalities in \eqref{eq:dudzeta:caseII} and using $\mathcal{J}^{+}=(\zeta^*,1)$ and $\mathcal{J}^{-} = (0,\zeta^*)$.

% \bibliographystyle{plain}
% \bibliography{references}

\begin{thebibliography}{10}

\bibitem{Barker2017}
T.~Barker, D.G.~Schaeffer, M.~Shearer, and J.M.N.T.~Gray.
\newblock Well-posed continuum equations for granular flow with compressibility and $\mu(i)$-rheology.
\newblock {\em Proc.~R.~Soc.~A Math.~Phys.~Eng.~Sci.}, 473(2201):20160846, 2017.

\bibitem{Bassi2020}
C.~Bassi, L.~Bonaventura, S.~Busto, and M.~Dumbser.
\newblock {A hyperbolic reformulation of the Serre-Green-Naghdi model for general bottom topographies}.
\newblock {\em Comput.~Fluids}, 212, 2020.

\bibitem{Bates2022}
P.D.~Bates.
\newblock Flood inundation prediction.
\newblock {\em Annu.~Rev.~Fluid Mech.}, 54(1):287--315, 2022.

\bibitem{bauerle2024rotationalinvariancehyperbolicityshallow}
M.~Bauerle, A.J.~Christlieb, M.~Ding, and J.~Huang.
\newblock On the rotational invariance and hyperbolicity of shallow water moment equations in two dimensions.
\newblock {\em SIAM J.~Math.~Anal.}, 57(1):1039--1085, 2025.

\bibitem{Burger2025b}
R.~B\"{u}rger, E.D.~Fernández-Nieto, J.~Garres-Díaz, and J.~Moya.
\newblock Well-balanced physics-based finite volume schemes for Saint-Venant–Exner-type models of sediment transport.
\newblock {\em Adv.~Water Resour.}, 206:105178, 2025.

\bibitem{Burger2025a}
R.~B\"{u}rger, E.D. Fernández-Nieto, and J.~Moya.
\newblock A multilayer shallow water model for tsunamis and coastal forest interaction.
\newblock {\em ESAIM: Math. Model. Numer. Anal.}, 59(2):1113--1144, 2025.

\bibitem{Cai2013}
Z.~Cai, Y.~Fan, and R.~Li.
\newblock Globally hyperbolic regularization of {G}rad's moment system in one dimensional space.
\newblock {\em Commun. Math. Sci.}, 11(2):547--571, 2013.

\bibitem{Canestrelli2009}
A.~Canestrelli, A.~Siviglia, M.~Dumbser, and E.~F. Toro.
\newblock Well-balanced high-order centred schemes for non-conservative hyperbolic systems.
Applications to shallow water equations with fixed and mobile bed.
\newblock {\em Adv. Water Resour.}, 32(6):834--844, 2009.

\bibitem{Careaga2024}
J.~Careaga and V.~Osores.
\newblock A multilayer shallow water model for polydisperse reactive sedimentation.
\newblock {\em Appl.~Math.~Model.}, 134:570--590, 2024.

\bibitem{castro2020well}
M.~J. Castro and C.~Par{\'e}s.
\newblock Well-balanced high-order finite volume methods for systems of balance laws.
\newblock {\em J.~Sci.~Comput.}, 82(2):1--48, 2020.

\bibitem{Ciallella2025}
M.~Ciallella, L.~Micalizzi, V.~Michel-Dansac, P.~\"{O}ffner, and D.~Torlo.
\newblock A high-order, fully well-balanced, unconditionally positivity-preserving finite volume framework for flood simulations.
\newblock {\em GEM, Int.~J.~Geomath.}, 16(1), 2025.

\bibitem{Cleary2017}
P.W.~Cleary, J.E.~Hilton, and M.D.~Sinnott.
\newblock Modelling of industrial particle and multiphase flows.
\newblock {\em Powder Technol.}, 314:232--252, 2017.

\bibitem{Coulomb1776}
C.A.~Coulomb.
\newblock Essai sur une application des regles des maximis et minimis a quelquels problemesde statique relatifs, a la architecture.
\newblock {\em Mem.~Acad.~Roy.~Div.~Sav.}, 7:343--387, 1776.

\bibitem{Fei2020}
J.B.~Fei, Y.X.~Jie, D.B.~Zhao, and B.Y.~Zhang.
\newblock Simulation of natural shallow avalanches with the $\mu(i)$ rheology.
\newblock {\em Bull.~Eng.~Geol.~Environ.}, 79(8):4123--4134, 2020.

\bibitem{FernandezNieto2016}
E.D.~Fern{\'{a}}ndez-Nieto, J.~Garres-D{\'{i}}az, A.~Mangeney, and G.~Narbona-Reina.
\newblock {A multilayer shallow model for dry granular flows with the $\mu$(I)-rheology: Application to granular collapse on erodible beds}.
\newblock {\em J.~Fluid Mech.}, 798:643--681, 2016.

\bibitem{FernandezNieto2018}
E.D.~Fernández-Nieto, J.~Garres-Díaz, A.~Mangeney, and G.~Narbona-Reina.
\newblock {2D} granular flows with the $\mu(i)$ rheology and side walls friction: {A} well-balanced multilayer discretization.
\newblock {\em J. Comput. Phys.}, 356:192--219, 2018.

\bibitem{Forterre2008}
Y.~Forterre and O.~Pouliquen.
\newblock Flows of dense granular media.
\newblock {\em Annu.~Rev.~Fluid Mech.}, 40(1):1--24, 2008.

\bibitem{Garres-Diaz2021}
J.~Garres-D{\'{i}}az, E.~D. Fern{\'{a}}ndez-Nieto, A.~Mangeney, and T.~{Morales de Luna}.
\newblock {A Weakly Non-hydrostatic Shallow Model for Dry Granular Flows}.
\newblock {\em J.~Sci.~Comput.}, 86(2):1--32, 2021.

\bibitem{Garres2020}
J.~Garres-Díaz, T.~Morales de~Luna, M.~J. Castro, and J.~Koellermeier.
\newblock Shallow water moment models for bedload transport problems.
\newblock {\em Commun. Comput. Phys.}, 11(3):435--467, 2021.

\bibitem{Garres2023}
J.~Garres-Díaz, C.~Escalante, T.~Morales de~Luna, and M.J.~Castro Díaz.
\newblock A general vertical decomposition of euler equations: Multilayer-moment models.
\newblock {\em Appl.~Numer.~Math.}, 183:236--262, 2023.

\bibitem{Jop2005}
P.~Jop, Y.~Forterre, and O.~Pouliquen.
\newblock Crucial role of sidewalls in granular surface flows: consequences for the rheology.
\newblock {\em J.~Fluid Mech.}, 541:167--192, 2005.

\bibitem{Jop2006}
P.~Jop, Y.~Forterre, and O.~Pouliquen.
\newblock A constitutive law for dense granular flows.
\newblock {\em Nature}, 441(7094):727--730, 2006.

\bibitem{Kelfoun2005}
K.~Kelfoun and T.H.~Druitt.
\newblock Numerical modeling of the emplacement of socompa rock avalanche, Chile.
\newblock {\em J. Geophys. Res. Solid Earth}, 110(B12), 2005.

\bibitem{Kern2010}
M.A.~Kern, P.A.~Bartelt, and B.~Sovilla.
\newblock Velocity profile inversion in dense avalanche flow.
\newblock {\em Ann.~Glaciol.}, 51(54):27--31, 2010.

\bibitem{Kitsikoudis2020}
V.~Kitsikoudis, B.P.J.~Becker, Y.~Huismans, P.~Archambeau, S.~Erpicum, M.~Pirotton, and B.~Dewals.
\newblock Discrepancies in flood modelling approaches in transboundary river systems: Legacy of the past or well-grounded choices?
\newblock {\em Water Resour.~Manage.}, 34(11):3465--3478, 2020.

\bibitem{Koellermeier2017b}
J.~Koellermeier.
\newblock {\em Derivation and numerical solution of hyperbolic moment equations for rarefied gas flows}.
\newblock dissertation, RWTH Aachen University, Aachen, 2017.

\bibitem{KoellermeierQian2022}
J.~Koellermeier and Q.~Huang.
\newblock Equilibrium stability analysis of hyperbolic shallow water moment equations.
\newblock {\em Math.~Method.~Appl.~Sci.}, 2022.

\bibitem{Pimentel2020}
J.~Koellermeier and E.~Pimentel.
\newblock Steady states and well-balanced schemes for shallow water moment equations with topography.
\newblock {\em Appl.~Math.~Comput.}, 427, 2022.

\bibitem{Koellermeier2020}
J.~Koellermeier and M.~Rominger.
\newblock Analysis and numerical simulation of hyperbolic shallow water moment equations.
\newblock {\em Commun.~Comput.~Phys.}, 28((3)):1038--1084, 2020.

\bibitem{Koellermeier2014}
J.~Koellermeier, R.~Pascal Schaerer, and M.~Torrilhon.
\newblock A framework for hyperbolic approximation of kinetic equations using quadrature-based projection methods.
\newblock {\em Kinet.~Relat.~Models}, 7(3):531--549, 2014.

\bibitem{Koellermeier2017}
J.~Koellermeier and M.~Torrilhon.
\newblock Numerical study of partially conservative moment equations in kinetic theory.
\newblock {\em Commun.~Comput.~Phys.}, 21(04)(4):981--1011, 2017.

\bibitem{kowalski2018moment}
J.~Kowalski and M.~Torrilhon.
\newblock Moment approximations and model cascades for shallow flow.
\newblock {\em Commun.~Comput.~Phys.}, 25, 2019.

\bibitem{Mangeney2003}
A.~Mangeney‐Castelnau, J.‐P.~Vilotte, M.O.~Bristeau, B.~Perthame, F.~Bouchut, C.~Simeoni, and S.~Yerneni.
\newblock Numerical modeling of avalanches based on Saint Venant equations using a kinetic scheme.
\newblock {\em J.~Geophys.~Res.~Solid Earth}, 108(B11):2527--2544, 2003.

\bibitem{Manning1891}
R.~Manning.
\newblock On the flow of water in open channels and pipes.
\newblock {\em Trans.~Instn.~Civ.~Eng.~Irel.}, 20:161--207, 1891.

\bibitem{Meister2016}
A.~Meister and S.~Ortleb.
\newblock A positivity preserving and well-balanced dg scheme using finite volume subcells in almost dry regions.
\newblock {\em Appl.~Math.~Comput.}, 272:259--273, 2016.

\bibitem{Neglia2021}
F.~Neglia, R.~Sulpizio, F.~Dioguardi, L.~Capra, and D.~Sarocchi.
\newblock Shallow-water models for volcanic granular flows: A review of strengths and weaknesses of titan2d and flo2d numerical codes.
\newblock {\em J.~Volcanol.~Geotherm.~Res.}, 410:107146, 2021.

\bibitem{pares2006numerical}
C.~Par{\'e}s.
\newblock Numerical methods for nonconservative hyperbolic systems: A theoretical framework.
\newblock {\em SIAM J.~Numer.~Anal.}, 44(1):300--321, 2006.

\bibitem{Sanvitale2016}
N.~Sanvitale and E.T.~Bowman.
\newblock Using piv to measure granular temperature in saturated unsteady polydisperse granular flows.
\newblock {\em Granul.~Matter}, 18(3), 2016.

\bibitem{Savage1989}
S.B.~Savage and K.~Hutter.
\newblock The motion of a finite mass of granular material down a rough incline.
\newblock {\em J.~Fluid Mech.}, 199:177--215, 1989.

\bibitem{Savage1991}
S.~B. Savage and K.~Hutter.
\newblock The dynamics of avalanches of granular materials from initiation to runout. Part I: Analysis.
\newblock {\em Acta Mech.}, 86(1--4):201--223, 1991.

\bibitem{Scholz2023}
U.~Scholz, J.~Kowalski, and M.~Torrilhon.
\newblock Dispersion in shallow moment equations.
\newblock {\em Commun.~Appl.~Math.~Comput.}, 6(4):2155--2195, 2023.

\bibitem{Verbiest2025}
R.~Verbiest and J.~Koellermeier.
\newblock Capturing vertical information in radially symmetric flow using hyperbolic shallow water moment equations.
\newblock {\em Commun.~Comput.~Phys.}, 37(3):810--848, 2025.

\bibitem{Yong1999}
W.-A.~Yong.
\newblock Singular perturbations of first-order hyperbolic systems with stiff source terms.
\newblock {\em J.~Differ.~Equ.}, 155(1):89--132, 1999.

\bibitem{Zhuang2025}
Y.~Zhuang, B.W.~McArdell, and P.~Bartelt.
\newblock Comparative analysis of $\mu(I)$ and Voellmy-type grain flow rheologies in geophysical mass flows: Insights from theoretical and real case
  studies.
\newblock {\em Nat.~Hazards Earth Syst.~Sci.}, 25(6):1901--1912, 2025.

\end{thebibliography}

\end{document}